\documentclass[a4paper, 11pt]{article}
\usepackage[margin=1.2in]{geometry}
\pdfoutput=1
\usepackage{soul}
\usepackage{mathtools, subfigure}

\usepackage[utf8]{inputenc}

\usepackage{amsfonts, amssymb,  bm, amsthm}
\usepackage{graphicx, color, epstopdf}
 \usepackage{caption}
\usepackage[colorlinks, 
            linkcolor=red, 
            anchorcolor=blue, 
            citecolor=green
            ]{hyperref}
\usepackage[inline]{showlabels}
\usepackage[title]{appendix}
\allowdisplaybreaks

\newtheorem{mylemma}{Lemma}[section]

\newtheorem{myremark}{Remark}[section]

\def\Xint#1{\mathchoice
  {\XXint\displaystyle\textstyle{#1}}%
  {\XXint\textstyle\scriptstyle{#1}}%
  {\XXint\scriptstyle\scriptscriptstyle{#1}}%
  {\XXint\scriptscriptstyle\scriptscriptstyle{#1}}%
  \!\int}
\def\XXint#1#2#3{{\setbox0=\hbox{$#1{#2#3}{\int}$}
    \vcenter{\hbox{$#2#3$}}\kern-.5\wd0}}

\def\bi{{\bf i}}

\def\XXint#1#2#3{{\setbox0=\hbox{$#1{#2#3}{\int}$}	\vcenter{\hbox{$#2#3$}}\kern-.5\wd0}}

\def\dashint{\Xint-}

\begin{document}
\title{An FFT-accelerated PML-BIE Solver for Three-Dimensional Acoustic Wave Scattering in Layered Media} 
\author{Hangya Wang$^1$ and Wangtao Lu$^1$}
\captionsetup[figure]{labelfont={bf},labelformat={default},labelsep=period,name={Figure}}
\maketitle

\footnotetext[1]{School of Mathematical Sciences, Zhejiang University, Hangzhou
  310027, China. Email: wangtaolu@zju.edu.cn, wangmy123@zju.edu.cn.  }
  
\begin{abstract}
This paper is concerned with three-dimensional acoustic wave scattering in two-layer media, where the two homogeneous layers are separated by a locally perturbed plane featuring an axially symmetric perturbation. A fast novel boundary integral equation (BIE) method is proposed to solve the scattering problem within a cylindrical perfectly matched layer (PML) truncation. We use PML-transformed Green’s functions to derive BIEs in terms of single- and double-layer potentials for the wave field and its normal derivative on the boundary of each truncated homogeneous region. These BIEs, combined with interface and PML boundary conditions, form a complete system that accurately approximates the scattering problem. An FFT-based approach is introduced to efficiently and accurately discretize the surface integral operators in the BIEs, where a new kernel splitting technique is developed to resolve instabilities arising from the complex arguments in Green’s functions. Numerical experiments demonstrate the efficiency and accuracy of the proposed method, as well as the exponential decay of truncation errors introduced by the PML.
\end{abstract}

\section{Introduction}
Wave scattering phenomena in layered media play a crucial role in various domains of physics and engineering, including acoustics, electromagnetics, and elastodynamics. These phenomena are fundamental to understanding wave propagation mechanisms over complex surfaces and near-field optical interactions, with critical applications ranging from geophysical exploration to nano-optics\cite{chew1999waves}. Recently, there has been increasing interest in axisymmetric structures\cite{helsing2017, ONEIL2018263}.
Despite the significance of these problems, various challenges limit the availability of accurate and efficient numerical solutions. 
Traditional numerical methods, such as the Finite Element Method (FEM) \cite{monk2003finite}, while widely implemented, become too expensive when the scatterer size significantly exceeds the wavelength. 
In contrast, the BIE method  \cite{colton2013} provides an effective alternative for various scattering scenarios \cite{bruno2012regularized,bruno2021,zhang2021accurate}. This approach has two main advantages: it reduces the problem dimension by one and inherently satisfies the outgoing radiation condition. For large-scale problems, the BIE method outperforms classical discretization methods, enabling faster and higher-order algorithms.

BIE solvers for layered medium scattering problems can be broadly categorized into two principal approaches.  The first approach uses the background Green's functions of the layered medium  \cite{Sommerfeld1909, Paulus2000,perez2014}, with integral equations formulated at local interfaces or boundaries. These Green's functions are defined only on the bounded region of the scattering surface, thereby eliminating the need for additional truncation. However, this approach presents a significant computational challenge due to the need to evaluate complex Sommerfeld integrals arising from the layered-medium Green's function and its derivatives.
A comprehensive computational cost analysis for evaluating the Sommerfeld integrals is available in  \cite{cai2002}.
The second approach employs free-space Green's functions and requires the formulation of integral equations on unbounded interfaces between different background medium layers. While this approach avoids intricate Sommerfeld integrals, special treatments are required to truncate the unbounded surface. These treatments include the
approximate truncation method \cite{meier2001,saillard2011}, the taper function method  \cite{miret2013,spiga2008,zhao2005}, and the windowed function method  \cite{bruno2016,monro2008,bruno2021}. The windowed Green function method \cite{bruno2016} attains superior performance through a correction and ensures uniformly fast convergence across all incident angles as the window size increases.

Recently,  \cite{lulu2018} proposed a BIE method based on PMLs for two-dimensional (2D) scattering problems in layered media. Similar to BIE methods based on the free-space Green’s function, this approach avoids evaluating expensive Sommerfeld integrals by formulating integral equations along the interfaces of the background layered medium. This method employs the PML-transformed free-space Green’s function to attenuate outgoing waves, enabling direct truncation of the unbounded scattering domain and leading to solutions that converge exponentially \cite{chen2010}. So far, this approach has been successfully extended to more complex structures, including step-scattering surfaces  \cite{lu2021}, anisotropic media  \cite{gao2022}, localized regenerative periodic structures  \cite{yu2022}. and 2D layered media with local smooth defects. Notably, \cite{luxu2023} extends this approach to three dimensions and develops efficient treatment for the hyper-singular boundary integral operators.

In this paper, we develop a novel FFT-accelerated PML-BIE method for 3D acoustic scattering problems in axially symmetric, piecewise homogeneous media. The computational domain comprises two homogeneous regions separated by a flat interface containing a finite, axially symmetric perturbation. Employing PML technology efficiently handles the unbounded domain.
Due to the axisymmetry, our method has significant advantages in the following aspects. Our approach is based on a direct boundary integral formulation derived from Green's representation theorem, which avoids the discretization of hyper-singular boundary integral operators. With axisymmetry, we apply Fourier transformation to reduce the surface integral equations into a series of decoupled curve integral equations \cite{wanglu2024}. The decoupling enables the use of the FFT to enhance computational efficiency significantly.
To accurately evaluate both regular and singular integrals, we employ high-order quadrature rules developed by Alpert\cite{Alpert1999}.
To address instability issues and obtain accurate PML-transformed modal Green’s functions, we propose a modified kernel-splitting technique based on function substitutions and numerical techniques.
A carefully designed graded mesh scheme  \cite{colton2013} that effectively captures solution behavior near singularities is employed to solve the geometric singularity. 
Numerical experiments demonstrate that the PML truncation of unbounded interfaces is effective, yielding high-accuracy results when appropriate PML parameters are selected.

The rest of this paper is organized as follows: Section 2 presents the formulation of scattering problems in piecewise constant penetrable media. Section 3 introduces boundary integral equations for the truncated PML problem and applies Fourier transformation to these equations. Section 4 details the numerical discretization scheme for the integral equations and several significant numerical techniques. Finally, Section 5 demonstrates the performance of our method through numerical experiments. 

\section{Problem formulation}
\begin{figure}[htb]
    \centering
    \includegraphics[width=0.6\textwidth]{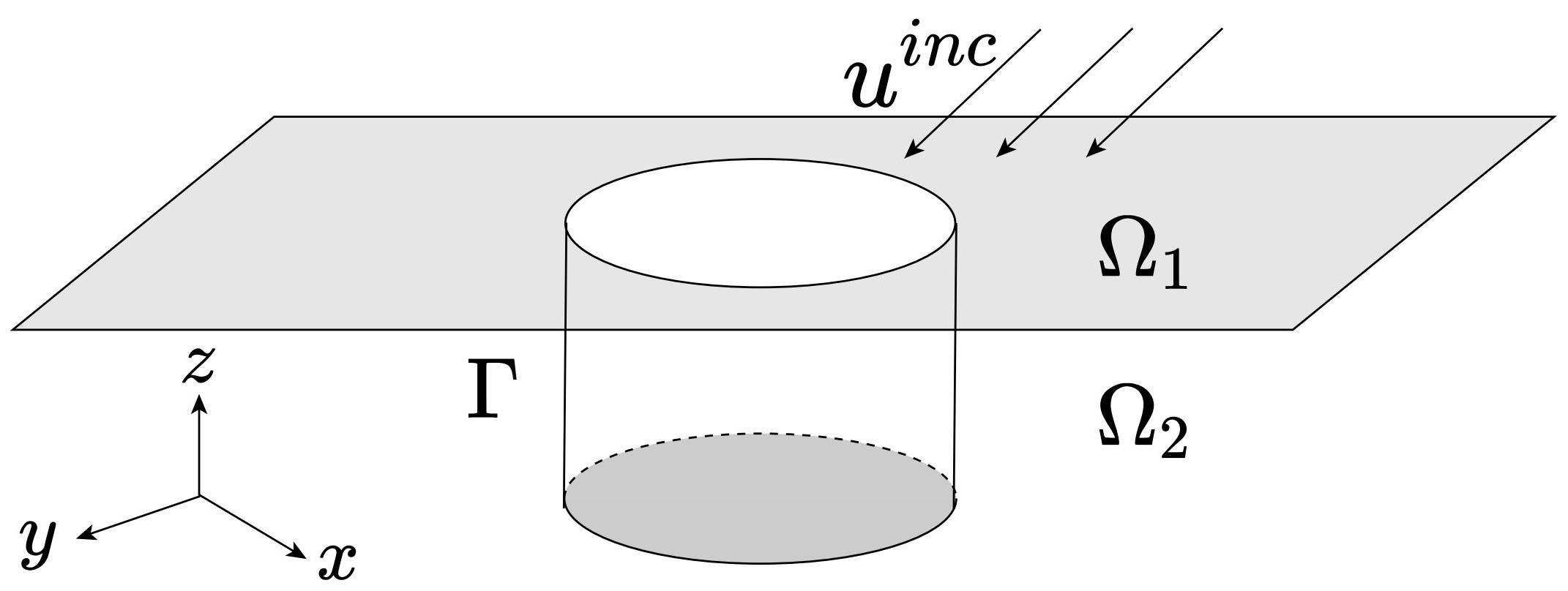}

    \caption{Profile of a 3D layered medium.} 
    \label{fg1}
\end{figure}

As illustrated in Figure \ref{fg1}, the layered medium consists of two homogeneous domains, denoted as $\Omega_1$ and $\Omega_2$. 
The interface $\Gamma$ separating $\Omega_1$ and $\Omega_2$ consists of an infinite flat surface $z=0$ and an axisymmetric perturbation surface. In what follows, we assume that the surface $\Gamma$ is piecewise smooth and contains a finite number of corners/edges. 

Let $u^{inc}$ represent a time-harmonic incident wave field with time-dependence $e^{-\bi \omega t}$. The total wave field after scattering is denoted by $u^{tot}$. In the standard Cartesian coordinate system, the governing 3D Helmholtz equation for $u^{tot}$ is
\begin{align}
\label{eq:helm}
    \Delta u^{tot} + \frac{\omega^2}{c^2(\boldsymbol{r})} u^{tot} =& 0, \quad{\rm in}\quad \mathbb{R}^3,
\end{align}
where $\boldsymbol{r}=(x,y,z)$, $\omega$ denotes the angular frequency and $\Delta = \partial^2_x+ \partial^2_y+ \partial^2_z$ represents the 3D Laplacian operator.
Here, $c(\boldsymbol{r})$ represents the speed of sound in the layered media, which is given by 
\[
c({\bm r}) = c_j, \quad \text{in} \quad \Omega_j ,\quad j=1,2.
\]
For simplicity, the density is assumed to be constant throughout the entire space. Let $\omega_j = \omega/c_j$, with $c_1=1$. Across the interface $\Gamma$, the total wave field $u^{tot}$ satisfies the continuity condition 
\begin{equation}
\label{eqbdc}
[u^{tot}] = 0, \quad {\rm and} \quad \left[\frac{\partial u^{tot}}{\partial \nu}\right]=0,
\end{equation}
where $[\cdot]$ denotes the jump of the quantity, and $\nu$ is the unit outer normal vector on $\Gamma$ towards $\Omega_2$. 
Since the total field $u^{tot}$ does not satisfy the outgoing radiation condition, we express it as the sum of
a reference wave field $u^{tot,0}$, defined as follows, and a scattered wave field $u^s=u^{tot}-u^{tot,0}$. 

In this paper, we consider two types of incident waves, a plane wave and a spherical wave excited by a point source. For plane wave incidence, we assume the incident wave is given by 
\begin{equation}
    \label{eqpw1}
    u^{inc}=e^{\bi \omega_1(-x \sin \phi^{inc} \cos \theta^{inc} -y \sin \phi^{inc} \sin \theta^{inc} -z \cos\phi^{inc})},
\end{equation}
where the polar angle $\phi^{inc} \in [0,\pi/2)$ and the azimuthal angle $\theta^{inc} \in [0,2\pi)$ .
Due to the axial symmetry of the layered structure, we set $\theta^{inc}=0$. This simplifies equation \eqref{eqpw1} to
\begin{equation}\label{eqpw}
    u^{inc}=e^{\bi \omega_1(-x \sin \phi^{inc} -z \cos\phi^{inc})}.
\end{equation}
The reference wave field $u^{tot,0}$ is defined as the solution of the incident wave $u^{inc}$ scattered by a flat interface at $z=0$. The explicit expression for $u^{tot,0}$ is given by
\begin{equation}
    u^{tot,0}=\left\{
    \begin{aligned}
        &e^{\bi \omega_1(-x \sin \phi^{inc} -z \cos\phi^{inc})}+R e^{\bi \omega_1(-x \sin \phi^{inc} +z \cos\phi^{inc})}, \quad &in \quad \Omega_1,\\
        &(R+1)e^{-\bi \omega_1 x \sin \phi^{inc} -\bi k^* z}, \quad &in \quad \Omega_2,
    \end{aligned}
    \right.
\end{equation}
where
\begin{align}
    R=\dfrac{2}{1+\frac{k^*}{\omega_1 \cos\phi^{inc}}}-1.
\end{align}
Here, $k^*$ is defined as
\begin{align}
     k^*= \left\{\begin{array}{ll}
\sqrt{\omega_2^2-\omega_1^2 \sin^2 \phi^{inc}}, & \omega_2 > |\omega_1 \sin \phi^{inc}|,\\
\bi\sqrt{\omega_1^2 \sin^2 \phi^{inc}-\omega_2^2}, & \omega_2 \leq |\omega_1 \sin \phi^{inc}|. \\
\end{array}\right.
\end{align}

For spherical wave incidence, equation \eqref{eq:helm} is modified to
\begin{align}
\label{eq:helm1}
    \Delta u^{tot} + \omega_j^2 u^{tot} =& -\delta(\boldsymbol{r},\boldsymbol{r^{*}}), \quad{\rm in}\quad \Omega_j, \quad j=1,2,
\end{align}
where $\boldsymbol{r^{*}}=(x^{*},y^{*},z^{*})$ represents the coordinates of the source point. The incident wave is expressed  as 
\begin{equation}
    u^{inc}=\frac{e^{\bi \omega_1|\boldsymbol{r}-\boldsymbol{r}^{*}|}}{4 \pi\left|\boldsymbol{r}-\boldsymbol{r}^{*}\right|}.
\end{equation}
For $\boldsymbol{r}^{*} \in \Omega_1$, the reference wave field $u^{tot,0}$ is given by
\begin{equation}
    u^{tot,0}=\left\{
    \begin{aligned}
        &u^{inc}, \quad & \text{in} \quad \Omega_1,\\
        &0, \quad &\text{in} \quad \Omega_2.
    \end{aligned}
    \right.
\end{equation}

Based on equation \eqref{eq:helm} and the definition of $u^{tot,0}$, the scattered field $u^s$ satisfies the Helmholtz equation
\begin{equation}
    \Delta u^s+\omega_j^2 u^s=0  , \quad{\rm in}\quad \Omega_j, \quad j=1,2.
\end{equation}
At infinity, the scattered field $u^s$ satisfies the following half-space Sommerfeld radiation condition
\begin{equation}
\lim _{|\boldsymbol{r}| \rightarrow \infty} |\boldsymbol{r}|\left(\partial_{|\boldsymbol{r}|}- \bi \omega_j \right) u^{s}=0, \quad{\rm in}\quad \Omega_j, \quad j=1,2,
\end{equation}
where $|\cdot|$ denote the Euclidean norm; see the equivalent upward propagating radiation condition \cite{Arens05,chandler98} for more details. Instead of computing $u^{tot}$  directly, we focus on computing the scattered field $u^{s}$, which satisfies the following boundary conditions on the interface $\Gamma$:
\begin{align}
\label{eqs1}
\left.u^{s,1}\right|_{\Gamma}-\left.u^{s,1}\right|_{\Gamma} & =-\left[u^{tot,0}\right] ,\\
\label{eqs2}
\left. \frac{\partial u^{s,1}}{\partial \nu}\right|_{\Gamma}-\left.\frac{\partial u^{s,2}}{\partial \nu}\right|_{\Gamma} & =-\left[ \frac{\partial u^{tot,0}}{\partial \nu}\right],
\end{align}
where $u^{s,j}$ denotes $u^s$ in $\Omega_j$, for $j = 1, 2$.
In the following section, we will introduce a PML-based BIE formulation for solving the scattering problems.

\section{The PML truncation and boundary integral equations.} 
This section introduces the PML truncation and
the boundary integral equations associated with the truncated PML problems. 
Without loss of generality, we restrict our analysis to the upper homogeneous domain $\Omega_1$. Henceforth, we denote $u^{s,1}$ and $\omega_1$ simply as $u^s$ and $\omega$, respectively.

\subsection{Boundary integral equations for the PML problem}
\begin{figure}
    \centering
    \subfigure[]{
    \includegraphics[width=0.45\linewidth]{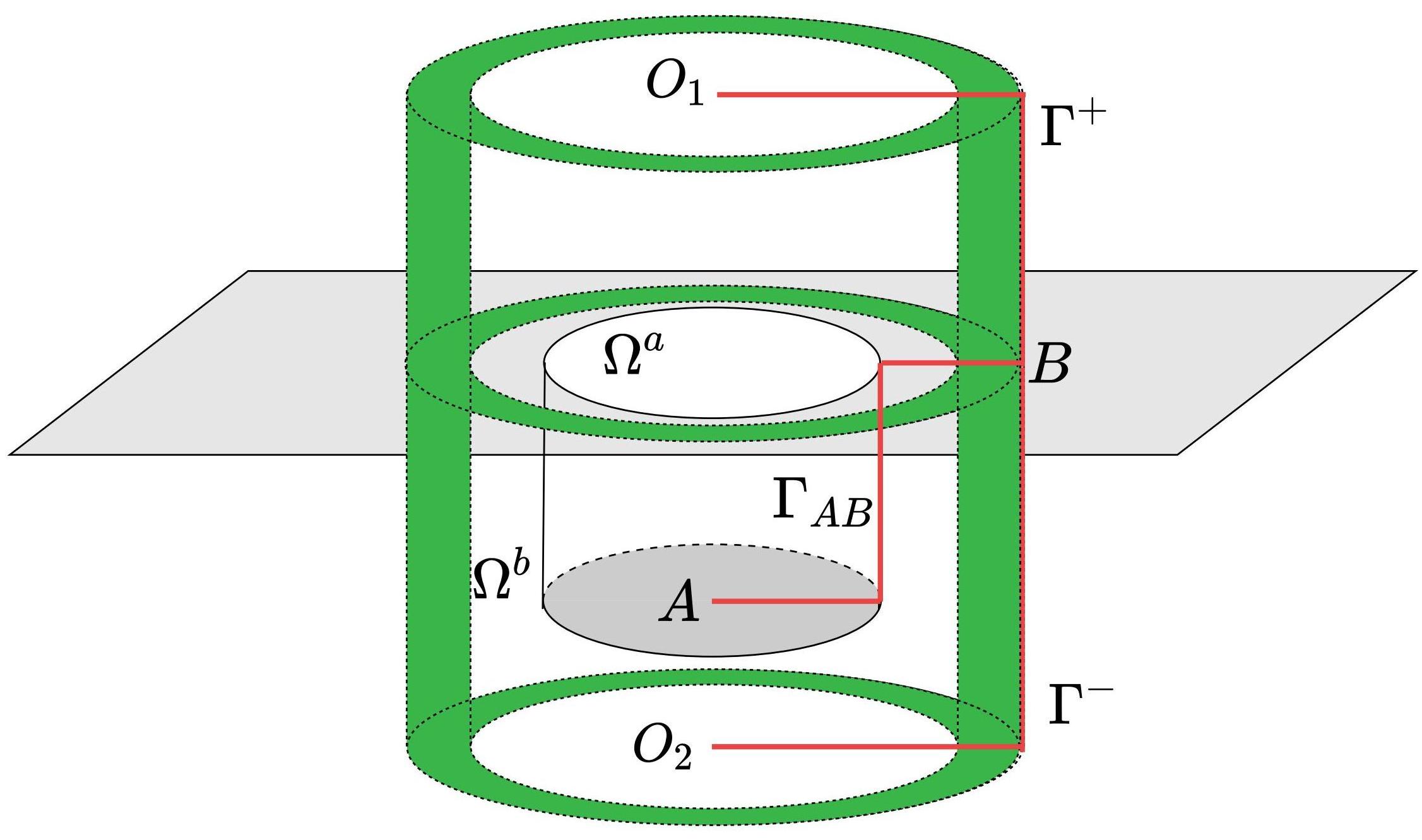}
    }
    \subfigure[]{
    \includegraphics[width=0.3\linewidth]{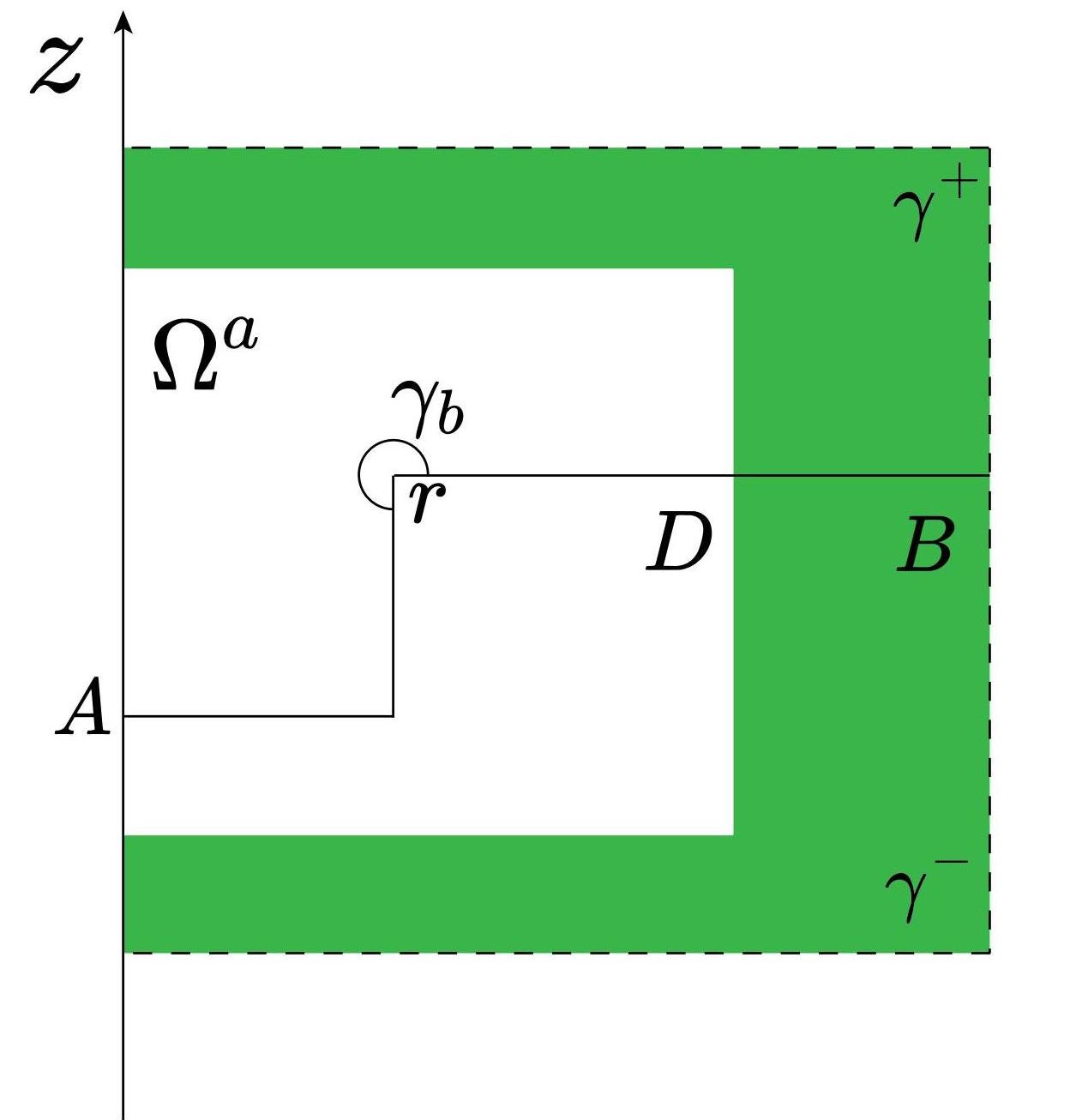}
    }
    \caption{ (a) The 3D illustration of the PML truncation. (b) The 2D cross-section of the PML truncation.}
    \label{fg2}
\end{figure}

As illustrated in Figure \ref{fg2}(a), a cylindrical PML, represented by the green region, is employed to enclose the perturbation and truncate the domain $\mathbb{R}^3$. 
Let $\Gamma_{AB}$, $\Gamma^+$ and $\Gamma^-$ denote the surfaces generated by rotating the red polylines connecting points $A$ to $B$, $B$ to $O_1$, and $B$ to $O_2$, respectively. The cylindrical domain consists of two bounded domains $\Omega^a$ and $\Omega^b$ with boundaries $\Gamma^a= \Gamma^+ \cup \Gamma_{AB}$ and $\Gamma^b=\Gamma^- \cup \Gamma_{AB}$, respectively.
Following \cite{Chew1994}, we introduce a complex coordinate stretching 
$\tilde{\boldsymbol{r}}(\boldsymbol{r})= 
(\tilde{\rho}\cos \theta,\tilde{\rho}\cos \theta,\tilde{z})$, which is defined as
\begin{align}
    \label{eqrhot}
    \tilde{{\rho}}({\rho})=\rho+\bi \int_0^\rho \sigma_1(t) dt, \\
    \label{eqzt}
    \tilde{{z}}({z})=z+\bi \int_0^z\sigma_2(t) dt,
\end{align}
where $\sigma_l(t)=\sigma_l(-t)$, $\sigma_l(t)=0$ for $|t|\leq a_l$ and $\sigma_l(t)>0$ for $|t|>a_l$, for $l=1,2$. 
In principle, $\sigma_l$ can be any positive function in the PML. For definiteness, taking $\sigma_1$ as an example, it is defined as
\begin{equation}
\label{eqsgm1}
\begin{aligned}
&\sigma_1\left(\rho\right)=\left\{\begin{array}{cc}
\frac{2 S f_1^p}{f_1^p+f_2^p}, & a_1 \leq \rho \leq a_1+T ,\\
S, & \rho>a_1+T, \\
0, &  a_1<\rho<a_1, \\
\sigma_1\left(-\rho\right), & \rho \leq-a_1,
\end{array}\right.\\
\end{aligned}
\end{equation}
where $p$ is a positive integer,
\begin{equation}
    f_1=\left(\frac{1}{2}-\frac{1}{p}\right) \bar{\rho}^3+\frac{\bar{\rho}}{p}+\frac{1}{2}, \quad f_2=1-f_1, \quad \bar{\rho}=\frac{\rho-\left(a_1+T\right)}{T}.
\end{equation}
The function $\sigma_1(\rho)$ maps the interval $[a_1,a_1+T]$ to $[0,S]$, where $T$ represents the PML thickness and $S$ determines the magnitude of $\sigma_1(\rho)$. Additionally, the derivatives of the function $\sigma_1$ vanish at $\rho = \pm a_1$ up to order $p-1$. A similar definition applies to $\sigma_2$.

According to Lemma 2.3 in \cite{Lassas2001}, the scattered field $u^s(\boldsymbol{r})$ can be analytically continued into the $\tilde{\Omega}^a=\{\tilde{\boldsymbol{r}}(\boldsymbol{r}) | \boldsymbol{r} \in  {\Omega}^a\}$. Within this domain, $u^s(\tilde{\boldsymbol{r}})$ satisfies the equation
\begin{equation}
\label{eqht}
\tilde{\Delta} u^s(\tilde{\boldsymbol{r}})+\omega^2 u^s(\tilde{\boldsymbol{r}})=0, \quad \text{in} \quad  \tilde{\Omega}^a,
\end{equation}
where $\tilde{\Delta}  =\dfrac{\partial^2 }{\partial \tilde{\rho}^2}+\dfrac{1}{\tilde{\rho}}\dfrac{\partial }{\partial \tilde{\rho}}+\dfrac{1}{\tilde{\rho}^2}\dfrac{\partial^2 }{\partial \theta^2}+\dfrac{\partial^2 }{\partial \tilde{z}^2}$. Applying Green's representation theorem, the solution of \eqref{eqht} can be expressed as
\begin{equation}
u^s(\tilde{\boldsymbol{r}})=\int_{\Gamma^a}\left\{G(\tilde{\boldsymbol{r}}, \boldsymbol{r}') \partial_\nu u^s(\boldsymbol{r}')-\partial_\nu G(\tilde{\boldsymbol{r}}, \boldsymbol{r}') u^s(\boldsymbol{r}')\right\} d \boldsymbol{r}'.
\end{equation}
Using the notation $\tilde{u}^s({\boldsymbol{r}})=u^s(\tilde{\boldsymbol{r}})$ and chain rule, equation \eqref{eqht} can be reorganized as 
\begin{equation}
\label{eqht2}
\nabla \cdot\left(\mathbf{A} \nabla \tilde{u}^s\right)+\omega^2 J \tilde{u}^s=0, \quad \text{on} \quad \Omega^a,
\end{equation}
where $\alpha_\rho\left(\rho\right)=1+i \sigma_1\left(\rho\right) $, $ \alpha_z\left(z\right)=1+i\sigma_2\left(z\right)$, $\mathbf{A}=\operatorname{diag}(\dfrac{\tilde{\rho}  \alpha_z }{\rho\alpha_\rho},\dfrac{\rho }{\tilde{\rho}}\alpha_\rho\alpha_z,\dfrac{\tilde{\rho}\alpha_\rho }{\rho\alpha_z})$, $J=\dfrac{{\tilde{\rho} }}{\rho}\alpha_\rho\alpha_z$.
According to the Lemma 2.8 in  \cite{Lassas2001}, the fundamental solution of \eqref{eqht2} is
\begin{equation}
    \tilde{G}(\boldsymbol{r}, \boldsymbol{r}^\prime)=G(\tilde{\boldsymbol{r}}, \tilde{\boldsymbol{r}}^\prime)=\frac{e^{\bi \omega|\tilde{\boldsymbol{r}}-\tilde{\boldsymbol{r}}^{\prime}|}}{4 \pi\left|\tilde{\boldsymbol{r}}-\tilde{\boldsymbol{r}}^{\prime}\right|}, \quad \tilde{\boldsymbol{r}} \neq \tilde{\boldsymbol{r}}^\prime,
\end{equation}
where the complexified distance function is defined as
\begin{equation}
    |\tilde{\boldsymbol{r}}-\tilde{\boldsymbol{r}}^\prime|=\sqrt{\tilde{\rho}^2+\tilde{\rho}^{\prime 2}-2 \tilde{\rho} \tilde{\rho}^{\prime} \cos \left(\theta-\theta^{\prime}\right)+\left(\tilde{z}-\tilde{z}^{\prime}\right)^2}.
\end{equation}
Throughout this paper, for the function $\sqrt{z}$, we choose the negative real axis as the branch cut, ensuring that the argument $\arg\sqrt{z}$ lies within $(-\pi/2,\pi/2]$.
Using the results in \cite{lulu2018}, we have the following integral representation
\begin{equation}\label{eqg1t}
    \tilde{u}^s(\boldsymbol{r})=\int_{\Gamma^a }\left\{\frac{\partial \tilde{u}^s}{\partial \nu_c}(\boldsymbol{r}^\prime) \tilde{G}(\boldsymbol{r} ,  \boldsymbol{r}^\prime)-\tilde{u}^s(\boldsymbol{r}^\prime) \frac{\partial \tilde{G}(\boldsymbol{r},  \boldsymbol{r}^\prime)}{\partial  \nu_c( \boldsymbol{r}^\prime)}\right\} d S(\boldsymbol{r}^\prime),  \quad \boldsymbol{r} \in \Omega^a ,
\end{equation}
where 
\begin{equation}
\nu_c=\mathbf{A}^T\nu= \left( \nu_\rho \cos\theta \dfrac{\tilde{\rho}\alpha_z}{\rho\alpha_\rho},
 \nu_\rho \sin\theta \dfrac{\tilde{\rho}\alpha_z}{\rho\alpha_\rho},
 \nu_z \dfrac{\tilde{\rho}\alpha_\rho}{\rho\alpha_z}\right),    
\end{equation}
and $\frac{\partial}{\partial \nu_c}=\nu_c \cdot \nabla$. 
For simplicity, $\nu_\rho \dfrac{\tilde{\rho}\alpha_z}{\rho\alpha_\rho}$ and $\nu_z \dfrac{\tilde{\rho}\alpha_\rho}{\rho\alpha_z}$ are denoted by $\nu_{c,\rho}$ and $\nu_{c,z}$, respectively.
As $\boldsymbol{r}$ approaches the boundary $\Gamma^a $, equation \eqref{eqg1t} becomes
\begin{equation}\label{eq4}
    \left(\tilde{\mathcal{K}}- \tilde{\mathcal{K}}^0 1\right) \tilde{u}^s=\tilde{\mathcal{S} }\partial_{\nu_c} \tilde{u}^s,\quad{\rm on}\quad \Gamma^a.
\end{equation}
where
\begin{align}
\label{eq:St}
    (\tilde{\mathcal{S}} \varphi)(\boldsymbol{r}) =&2 \int_{\Gamma^a } \tilde{G}(\boldsymbol{r},  \boldsymbol{r}^\prime) \varphi(\boldsymbol{r}^\prime) d S(\boldsymbol{r}^\prime),  \quad \boldsymbol{r} \in \Gamma^a ,\\ 
\label{eq:Kt}
    (\tilde{\mathcal{K}} \varphi)(\boldsymbol{r}) =&2 \dashint_{\Gamma^a } \frac{\partial \tilde{G}(\boldsymbol{r},  \boldsymbol{r}^\prime)}{\partial  \nu_c(\boldsymbol{r}^\prime)} \varphi(\boldsymbol{r}^\prime) d S(\boldsymbol{r}^\prime),  \quad \boldsymbol{r} \in \Gamma^a,\\
    \label{eq:K0t}
    (\tilde{\mathcal{K}^0} \varphi)(\boldsymbol{r}) =&2 \dashint_{\Gamma^a } \frac{\partial \tilde{G}_0(\boldsymbol{r},  \boldsymbol{r}^\prime)}{\partial  \nu_c(\boldsymbol{r}^\prime)} \varphi(\boldsymbol{r}^\prime) d S(\boldsymbol{r}^\prime),  \quad \boldsymbol{r} \in \Gamma^a.
\end{align}
Here, $\dashint$ denotes the Cauchy principal value integral and $\tilde{G}_0(\boldsymbol{r}, \boldsymbol{r}^\prime)=\frac{1}{4 \pi\left|\tilde{\boldsymbol{r}}-\tilde{\boldsymbol{r}}^{\prime}\right|}$ represents the Green's function of the complexified Laplace equation
\begin{equation}
    \nabla \cdot\left(\mathbf{A} \nabla \tilde{u}^s_{n,0}(\boldsymbol{r})\right)=0.
\end{equation}
Finally, $\tilde{\mathcal{K}}^0 1(\boldsymbol{r})$ is defined as
\begin{equation}
\label{eqk01_2t}
    \tilde{\mathcal{K}}^0 1({\bm r})=2 \dashint_{\partial \Omega^a } \frac{\partial \tilde{G}_0(\boldsymbol{r},  \boldsymbol{r}^\prime)}{\partial  \nu_c(\boldsymbol{r}^\prime)}  d S(\boldsymbol{r}^\prime),  \quad \boldsymbol{r} \in \Gamma^a.
\end{equation}

\subsection{Truncation of boundary integral equations onto $\Gamma_{AB}$}
The PML is an artificial absorbing region designed to effectively absorb outgoing waves, ensuring that the scattered field $u^s$ and its normal derivative $\partial_{\nu} u^s$ decay rapidly within the PML before reaching its outer boundary. Due to the absorption properties of the PML, it is appropriate to approximate $\tilde{u}^s \approx 0$ and $\partial_{\nu_c}\tilde{u}^s \approx 0$ on $\Gamma^+$. Equation \eqref{eq4}  can be truncated as
\begin{equation}
    \label{eqkbab}
    (\tilde{\mathcal{K}}_{AB} \tilde{u}^s  - {\tilde{\mathcal{K}}}^{0} 1 ) \tilde{u}^s = \tilde{\mathcal{S}}_{AB} 
\partial_{\nu_c} \tilde{u}^s , \quad \text{on} \quad \Gamma_{AB},
\end{equation}
where $\tilde{\mathcal{K}}_{AB}$ and $\tilde{\mathcal{S}}_{AB}$ are the same as $\tilde{\mathcal{K}}$ and $\tilde{\mathcal{S}}$, respectively, but with the integral domain replaced by $\Gamma_{AB}$. However, the truncation of $\tilde{\mathcal{K}}^{0} 1$ onto $\Gamma_{AB}$ is not straightforward.

To analyze the operator $\tilde{\mathcal{K}}^{0}1$, we first consider the operator 
\begin{equation}
\label{eqk01_2}
    \mathcal{K}^0 1({\bm r})=2 \dashint_{\partial \Omega^a } \frac{\partial G_0(\boldsymbol{r},  \boldsymbol{r}^\prime)}{\partial  \nu(\boldsymbol{r}^\prime)}  d S(\boldsymbol{r}^\prime),  \quad \boldsymbol{r} \in \Gamma^a ,
\end{equation}
which is defined as in \eqref{eqk01_2t} but with $\tilde{G}_0(\boldsymbol{r}, \boldsymbol{r}^\prime)$ replaced by ${G}_0(\boldsymbol{r}, \boldsymbol{r}^\prime)=\frac{1}{4 \pi\left|{\boldsymbol{r}}-{\boldsymbol{r}}^{\prime}\right|}$.
For the operator ${\mathcal{K}}^{0}1$, we have the following Lemma.
\begin{mylemma}
\label{lem:K01}
For any $\boldsymbol{r} \in \Gamma^a $, we have
\begin{equation}
    \label{eqk0_1}
    \mathcal{K}^{0}1={1- \cos \phi_0},
\end{equation}
where $\phi_0$ is half the angle of the curve $\gamma_b$, see $\gamma_b $ in {Figure \ref{fg2}(b)}.
\begin{proof}
Using Green’s identity, equation \eqref{eqk01_2} can be rewritten as 
\begin{align}
    \label{eqk01t}
    {\mathcal{K}}^0 1(\boldsymbol{r})&=\lim _{r_c \rightarrow 0^{+}} 2 \int_{\partial B(\boldsymbol{r}, r_c) \cap \overline{\Omega^a }} \partial_{\nu} G_0(\boldsymbol{r}, \boldsymbol{r}^\prime) d S(\boldsymbol{r}^\prime) \\
    \label{eqk03}
    &=\lim _{r_c \rightarrow 0^{+}} \int_{\partial B(\boldsymbol{r}, r_c) \cap \overline{\Omega^a }} \frac{\boldsymbol{\nu} \left(\boldsymbol{r}^{\prime}\right) \cdot\left(\boldsymbol{r}-\boldsymbol{r}^{\prime}\right)}{2\pi \left|\boldsymbol{r}-\boldsymbol{r}^{\prime}\right|^3}  \rho^\prime d S(\boldsymbol{r}^\prime),
\end{align}
where $\partial B (\boldsymbol{r}, r_c)$ is a sphere centered at $\boldsymbol{r}$ with radius $r_c$, and the unit normal vector
$\boldsymbol{\nu}$ points toward $\Omega^b $. 

Next, we parameterize $\partial B(x, r_c) $ by 
\begin{equation}
\label{eqk04}
    \boldsymbol{r}^\prime=\boldsymbol{r}+r_c(\sin\phi \cos\theta ,\sin\phi \sin\theta ,\cos\phi ),
\end{equation}
for $\theta \in\left[-\pi, \pi\right]$ and $\phi \in\left[0, \pi\right]$. Let $\gamma_b $ denotes the projection of ${\partial B(\boldsymbol{r}, r_c) \cap \overline{\Omega^a }}$ on $\rho$-$z$ half-plane.
With coordinate axis rotation and the parameterisation \eqref{eqk04}, equation \eqref{eqk03} becomes
\begin{align}
    {\mathcal{K}}^0 1(\boldsymbol{r})&=\frac{1}{{2 \pi}} \int_{\gamma_b} \int_{-\pi}^{\pi} \frac{r_c}{r_c^3}\rho_s^2 \sin\phi d \theta d\phi \\
&=\int_{0}^{\phi_0}\sin\phi d\phi\\
&={1-\cos\phi_0},
\end{align}
where $\rho_s$ denotes the radial distance in the spherical coordinate system, and \eqref{eqk0_1} follows.
\end{proof}
\end{mylemma}

For $\tilde{\mathcal{K}}^{0}1$, we have the following Lemma.
\begin{mylemma}
\label{lem:K02}
For any ${{\boldsymbol{r}}} \in \Gamma^a $, we have
\begin{equation}
    \label{eqk0_1t}
     \tilde{\mathcal{K}}^{0}1={1- \cos\phi_0}.
\end{equation}
where $\phi_0$ is half the angle of the curve $\gamma_b$.
\begin{proof}
Using Green’s identity, we can easily see that 
\begin{align}
    \label{eqk03t}
    {\tilde{\mathcal{K}}}^0 1({\boldsymbol{r}})&=\lim _{r_c \rightarrow 0^{+}}  \int_{\partial B({\boldsymbol{r}}, r_c) \cap \overline{\Omega^a }} \frac{\boldsymbol{{\nu_c}} \left(\tilde{\boldsymbol{r}}^{\prime}\right) \cdot\left(\tilde{\boldsymbol{r}}-\tilde{\boldsymbol{r}}^{\prime}\right)}{ 2\pi\left |\tilde{\boldsymbol{r}}-\tilde{\boldsymbol{r}}^{\prime}\right|^3}  \rho^\prime d s({\boldsymbol{r}}^\prime)
\end{align}
where $\partial B ({\boldsymbol{r}}, r_c)$ is a sphere centered at $\boldsymbol{r}$ with radius $r_c$.
Let the points $\boldsymbol{r}$ and $\boldsymbol{r}^\prime$ be denoted by $({r}_1,{r}_2,{r}_3)$ and $({r}_1^\prime,{r}_2^\prime,{r}_3^\prime)$, respectively.
Similarly, $\tilde{\boldsymbol{r}}$ and $\tilde{\boldsymbol{r}}^\prime$ are denoted by
$(\tilde{r}_1,\tilde{r}_2,\tilde{r}_3)$ and $(\tilde{r}_1^\prime,\tilde{r}_2^\prime,\tilde{r}_3^\prime)$, respectively.

By \eqref{eqrhot}, \eqref{eqzt} and \eqref{eqk04}, we have 
\begin{align}
    \tilde{r}_1^\prime-\tilde{r}_1 &=\int_{r_1}^{r^\prime_1} \alpha_\rho(\sqrt{s^2+r_2^2})  d s \\
    &=\int_{r_1}^{r_1+r_c \sin\phi \cos\theta} \alpha_\rho(\sqrt{s^2+r_2^2})ds\\ &=\alpha_\rho\left(\sqrt{r_1^2+r_2^2}\right) r_c \sin\phi \cos\theta+O\left(r_c^2\right).
\end{align}
Similarly,
\begin{align}
    \tilde{r}_2^\prime-\tilde{r}_2 &=\alpha_\rho\left(\sqrt{r_2^2+r_1^2}\right) r_c \sin\phi \sin\theta+O\left(r_c^2\right) \\
    \tilde{r}_3^\prime-\tilde{r}_3, &=\alpha_z\left(r_3\right)  r_c \cos\phi +O\left(r_c^2\right).
\end{align}

With coordinate axis rotation, as $r_c$ approaches $0$,  equation \eqref{eqk03t} becomes 
\begin{align}
    {\tilde{\mathcal{K}}}^0 1(\tilde{\boldsymbol{r}})
    &=\frac{1}{{2 \pi}} \int_{\gamma_b} \int_{-\pi}^{\pi} \frac{\alpha_\rho^2  \sin^2(\phi)+\alpha_z^2  \cos^2(\phi)} {\alpha_\rho^2  \sin^2(\phi)+\alpha_z^2  \cos^2(\phi)} \frac{\rho_s^2}{r_c^2} \sin\phi d \theta d\phi \\
&=\int_{0}^{\phi_0}\sin\phi d\phi\\
&={1-\cos\phi_0},
\end{align}
and \eqref{eqk0_1t} follows.
\end{proof}
\end{mylemma}

Comparing \eqref{eqk0_1} and \eqref{eqk0_1t} gives
\begin{align}
    \tilde{\mathcal{K}}^{0}1 (\boldsymbol{r})&={\mathcal{K}}^{0}1 (\boldsymbol{r})\\
    \label{k01bc}
    &={\mathcal{K}}_{AB}^{0}1 (\boldsymbol{r})+ {\mathcal{K}}_{+}^{0}1 (\boldsymbol{r}),
\end{align}
where the definitions of ${\mathcal{K}}_{AB}^{0}$ and ${\mathcal{K}}_{+}^{0}$ are the same as ${\mathcal{K}}^{0}$ but with the integral domain replaced by $\gamma_{AB}$ and $\gamma^{+}$, respectively. 

\begin{figure}
    \centering
    \includegraphics[width=0.5\linewidth]{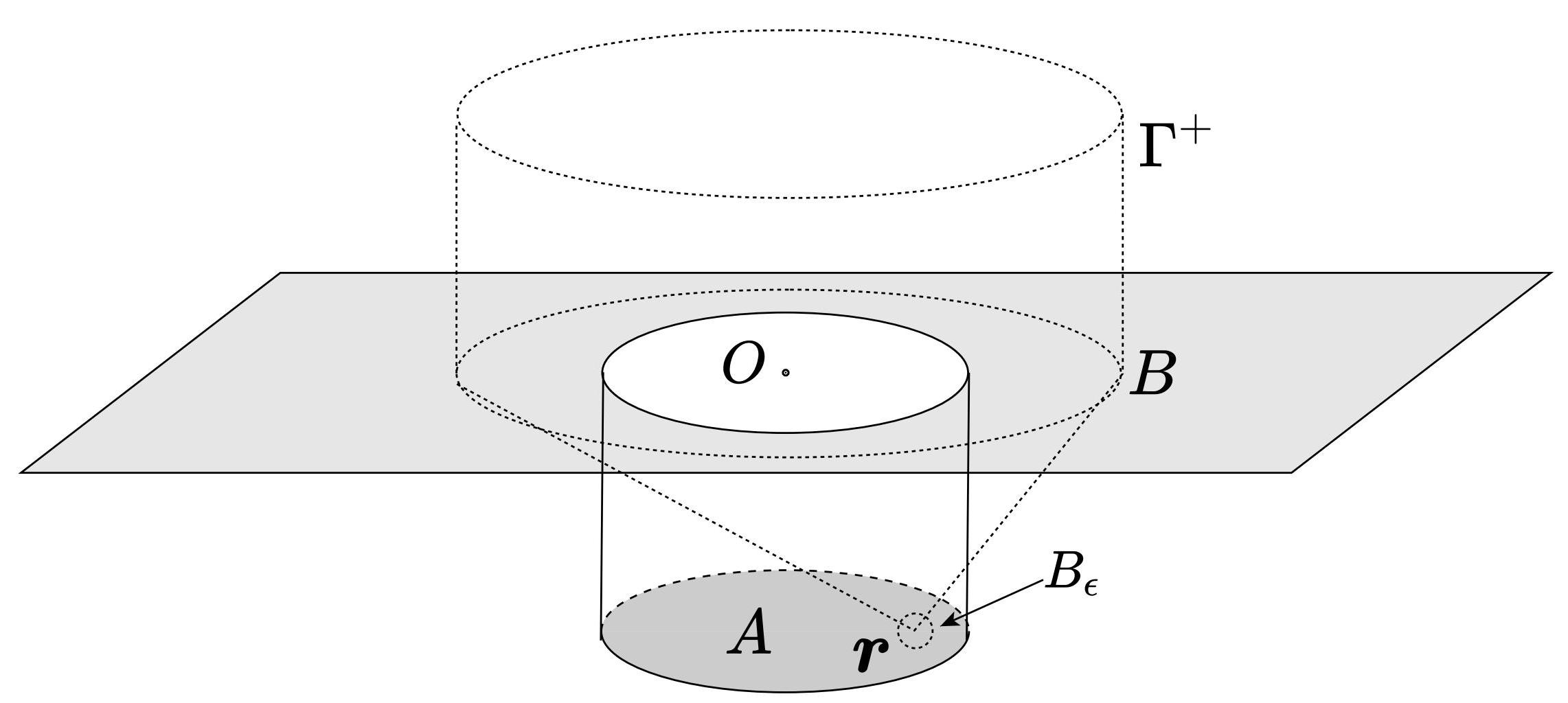}
    \caption{Geometric configuration for computing ${\mathcal{K}}_{+}^{0}1 (\boldsymbol{r})$}
    \label{fg8}
\end{figure}

To illustrate that the operator ${\mathcal{K}}_{+}^{0}1 (\boldsymbol{r})$ can be computed concisely, we give an auxiliary geometric configuration, as depicted in Figure \ref{fg8}. Let $\Gamma_{OB}$ denote the circle centered at origin $O$ with radius $l_{OB}$, where $l_{OB}$ is the length of $OB$. Let $\Gamma_{\boldsymbol{r}}$ denote the conical surface with $\Gamma_{OB}$ as the base and $\boldsymbol{r}$ as the vertex, where the point $\boldsymbol{r}=(\rho \cos \theta,\rho \sin \theta,z)$.  Let $\partial B_{\epsilon}$ represents the sphere centered at $\boldsymbol{r}$ with radius $\epsilon$. The intersection of the surface $\partial B_{\epsilon}$ and the conical surface $\Gamma_{\boldsymbol{r}}$ forms the surface $\Gamma_{\epsilon}$. In the following lemma, we provide a precise representation for ${\mathcal{K}}_{+}^{0}1 (\boldsymbol{r})$.

\begin{mylemma}
For any $\boldsymbol{r} $ on $ \Gamma_{AB} $, we have
\begin{align}
    \label{k01bc1}
    {\mathcal{K}}_{+}^{0}1 (\boldsymbol{r})
    &=\frac{1}{2 \pi \epsilon^2}  \int_{\Gamma_{\epsilon}}  \rho_s^2 \sin \phi d S(\boldsymbol{r}^\prime),
\end{align}
where $\rho_s$ denotes the radial distance in the spherical coordinate system.
\begin{proof}

According to the definition of ${\mathcal{K}}_{+}^{0}1$, we have
\begin{align}
    {\mathcal{K}}_{+}^{0}1 (\boldsymbol{r})
    \label{eqbdc1}
    &=\frac{1}{{2 \pi}} \int_{\Gamma^+} \frac{\boldsymbol{\nu} \left(\boldsymbol{r}^{\prime}\right) \cdot\left(\boldsymbol{r}-\boldsymbol{r}^{\prime}\right)}{\left|\boldsymbol{r}-\boldsymbol{r}^{\prime}\right|^3}  \rho^\prime d S(\boldsymbol{r}^\prime).
\end{align}
Applying the divergence theorem yields
\begin{equation}
    \begin{aligned}
        \int_{\Gamma^+} \frac{\boldsymbol{\nu} \left(\boldsymbol{r}^{\prime}\right) \cdot\left(\boldsymbol{r}-\boldsymbol{r}^{\prime}\right)}{\left|\boldsymbol{r}-\boldsymbol{r}^{\prime}\right|^3}  \rho^\prime d S(\boldsymbol{r}^\prime)+\int_{\Gamma_{\boldsymbol{r}}} \frac{\boldsymbol{\nu} \left(\boldsymbol{r}^{\prime}\right) \cdot\left(\boldsymbol{r}-\boldsymbol{r}^{\prime}\right)}{\left|\boldsymbol{r}-\boldsymbol{r}^{\prime}\right|^3}  \rho^\prime d S(\boldsymbol{r}^\prime) \\
    +\int_{\Gamma_{\epsilon}} \frac{\boldsymbol{\nu} \left(\boldsymbol{r}^{\prime}\right) \cdot\left(\boldsymbol{r}-\boldsymbol{r}^{\prime}\right)}{\left|\boldsymbol{r}-\boldsymbol{r}^{\prime}\right|^3}  \rho^\prime d S(\boldsymbol{r}^\prime)=0.
    \end{aligned}
\end{equation}
For ${\boldsymbol{r}^{\prime}} \in \Gamma_{\boldsymbol{r}}$, vectors $\boldsymbol{\nu}({\boldsymbol{r}^{\prime}})$ and $\boldsymbol{r}-\boldsymbol{r}^{\prime}$ are orthogonal, leading to
\begin{equation}
\label{eqbdc2}
    \int_{\Gamma_{\boldsymbol{r}}} \frac{\boldsymbol{\nu} \left(\boldsymbol{r}^{\prime}\right) \cdot\left(\boldsymbol{r}-\boldsymbol{r}^{\prime}\right)}{\left|\boldsymbol{r}-\boldsymbol{r}^{\prime}\right|^3}  \rho^\prime d S(\boldsymbol{r}^\prime)=0.
\end{equation}
Combining equations \eqref{eqbdc1} and \eqref{eqbdc2} yields
\begin{equation}
\label{eqk0j}
    {\mathcal{K}}_{+}^{0}1 (\boldsymbol{r})
    =\frac{-1}{{2 \pi}} \int_{\Gamma_{\epsilon}} \frac{\boldsymbol{\nu} \left(\boldsymbol{r}^{\prime}\right) \cdot\left(\boldsymbol{r}-\boldsymbol{r}^{\prime}\right)}{\left|\boldsymbol{r}-\boldsymbol{r}^{\prime}\right|^3}  \rho^\prime d S(\boldsymbol{r}^\prime).
\end{equation}
The surface $\Gamma_{\epsilon}$ is parameterized as
\begin{equation}
    \boldsymbol{r}^\prime=\boldsymbol{r}+\epsilon(\sin\phi \cos\theta ,\sin\phi \sin\theta ,\cos\phi ),
\end{equation}
so that $\boldsymbol{v}(\boldsymbol{r}^\prime)=(\sin\phi \cos\theta ,\sin\phi \sin\theta ,\cos\phi )$. Thus, equation \eqref{eqk0j} becomes
\begin{align}
    {\mathcal{K}}_{+}^{0}1 (\boldsymbol{r}) &=\frac{-1}{{2 \pi}}  \int_{\Gamma_{\epsilon}}\frac{\boldsymbol{v} \left(\boldsymbol{r}^\prime\right) \cdot\left(\boldsymbol{r}-\boldsymbol{r}^\prime\right)}{\left|\boldsymbol{r}-\boldsymbol{r}^\prime\right|^3}  \rho_s^2\sin \phi d S(\boldsymbol{r}^\prime)\\
    &=\frac{1}{2 \pi \epsilon^2}  \int_{\Gamma_{\epsilon}} \rho_s^2 \sin \phi d S(\boldsymbol{r}^\prime),
\end{align}
and this completes the proof.
\end{proof}
\end{mylemma}
\begin{myremark}
The integral $\int_{\Gamma_{\epsilon}} \rho_s^2 \sin \phi d S(\boldsymbol{r}^\prime)$ and $4 \pi \epsilon^2$ are the areas of surfaces 
$\Gamma_\epsilon$ and $\partial B_{\epsilon}$, respectively.
The areas of surface $\Gamma_\epsilon$ can be computed numerically via the Alpert quadrature rule. Hence, the computation of ${\mathcal{K}}_{+}^{0}1(\boldsymbol{r})$ is independent of $\Gamma^+$.   
\end{myremark}

By \eqref{k01bc} and \eqref{k01bc1}, equation \eqref{eqkbab}  can be truncated onto the interface $\Gamma_{AB}$ only, say
\begin{equation}
    \label{eq5}
    \tilde{\mathcal{K}}_{AB} \tilde{u}^s  - ({\mathcal{K}}^{0}_{AB}1 +{\mathcal{K}}_{+}^{0}1 ) \tilde{u}^s = \tilde{\mathcal{S}}_{AB} 
\partial_{\nu_c} \tilde{u}^s , \quad \text{on} \quad \Gamma_{AB}.
\end{equation}
Thus, we define the PML-transformed NtD map $\tilde{\mathcal{N} }= (\tilde{\mathcal{K}}_{AB}  - {\mathcal{K}}^{0}_{AB}1 -{\mathcal{K}}_{+}^{0}1 )^{-1}\tilde{\mathcal{S}}_{AB} $, mapping $\partial_{\nu_c} \tilde{u}^s $ to $ \tilde{u}^s $ on $\Gamma_{AB}$.
Numerically discretizing the integral operators in \eqref{eq5} approximates the NtD map $\tilde{\mathcal{N} }$ on $\Gamma_{AB}$.

\subsection{Fourier transformation of boundary integral equations }
Direct high-order discretizations of boundary integral equations are non-trivial and usually expensive for complex geometries in three dimensions. However, for axisymmetric domain $\Omega^a$,  
applying Fourier transformations converts surface integrals into a series of line integrals and yields an efficient algorithm. Thus, the first step in our discretization scheme is the Fourier transformation. 

The $n$-th azimuthal Fourier mode of generic functions $f(\boldsymbol{r})$ and $F(\boldsymbol{r},\boldsymbol{r}')$ can be expressed as  
\begin{equation}
\label{fmode1}
f_n(r)=\frac{1}{\sqrt{2\pi}}\int_{-\pi}^\pi e^{-\bi n \theta} f(\boldsymbol{r}) d \theta,
\end{equation}
and 
\begin{equation}
\label{fmode3}
F_n\left(r,  r^{\prime}\right)= \frac{1}{\sqrt{2\pi}}\int_{-\pi}^\pi e^{-\bi n(\theta -\theta^\prime) } F\left(\boldsymbol{r},  \boldsymbol{r}^{\prime}\right) d (\theta -\theta^\prime),    
\end{equation}
respectively, where $r=(\rho,z)$ and $r'=(\rho',z')$ are the projections of ${\bm r}$ and ${\bm r}'$ on the $\rho$-$z$ half-plane, respectively. 
Using the convolution formula,  the azimuthal Fourier modes of the product function $(fg)({\bm r})$ can be expressed in terms of $f_n(r)$ and $g_n(r)$, say 
\begin{equation}
    \label{eqcon}
    (f g)_n(r)=\frac{1}{\sqrt{2\pi}}\sum_{m=-\infty}^{\infty}f_m(r)g_{n-m}(r).
\end{equation}

According to the PML settings in Section 3.2, the complex transformation $\tilde{\boldsymbol{r}}(\boldsymbol{r})= 
(\tilde{\rho}\cos \theta,\tilde{\rho}\cos \theta,\tilde{z})$ is independent of $\theta$. Let $\gamma^a$ denotes the projections of $\Gamma^a$ on the $\rho$-$z$ half-plane. The kernels of the operators  $\tilde{\mathcal{S}}(\boldsymbol{r}, \boldsymbol{r}^\prime)$, $\tilde{\mathcal{K}} (\boldsymbol{r}, \boldsymbol{r}^\prime)$ and $\tilde{\mathcal{K}}^0 (\boldsymbol{r}, \boldsymbol{r}^\prime)$ are denoted by $\tilde{S}(\boldsymbol{r}, \boldsymbol{r}^\prime)$, $\tilde{K}(\boldsymbol{r}, \boldsymbol{r}^\prime)$ and $\tilde{K}^0(\boldsymbol{r}, \boldsymbol{r}^\prime)$, respectively. For any $n\in\mathbb{Z}$, integration of \eqref{eq4} yields

\begin{equation}
\label{eqntd1pml}
    \tilde{\mathcal{K}}_n \tilde{u}^s_n -\frac{1}{\sqrt{2 \pi}} [(\tilde{\mathcal{K}}^{0}1) \tilde{u}^s]_n= \tilde{\mathcal{S}}_n 
\partial_{\nu_c} \tilde{u}^s_n , \quad \text{on} \quad \gamma^a,
\end{equation}
where
\begin{align}
    \label{eqsnpml}
    (\tilde{\mathcal{S}}_n \varphi)(r) = 2\int_{\gamma^a } \tilde{S}_n(r,r')\varphi(r')\rho'ds(r'),\\
    \label{eqknpml}
    (\tilde{\mathcal{K}}_n \varphi)(r) = 2\dashint_{\gamma^a } \tilde{K}_n(r,r')\varphi(r')\rho'ds(r'),\\
    \label{eqk0pml}
    (\tilde{\mathcal{K}}_n^0 \varphi)(r) = 2\dashint_{\gamma^a } \tilde{K}^0_n(r,r')\varphi(r')\rho'ds(r'),
\end{align}
with $\tilde{S}_n(r,  r^{\prime})$, $\tilde{K}_n(r,  r^{\prime})$ and $\tilde{{K}}_n^0(r,  r^{\prime})$ denote the $n$-th azimuthal Fourier modes of $\tilde{S}(\boldsymbol{r}, \boldsymbol{r}^\prime)$, $\tilde{K}(\boldsymbol{r}, \boldsymbol{r}^\prime)$ and $\tilde{K}^0(\boldsymbol{r}, \boldsymbol{r}^\prime)$, respectively. Using the convolution formula and the transform
\begin{equation}
    (\tilde{\mathcal{K}}^0 1)_n(r)
   = \sqrt{2\pi}(\tilde{\mathcal{K}}^0_n 1_n) (r)=\delta_{n0} (\tilde{\mathcal{K}}^0 1)_0(r),
\end{equation}
equation \eqref{eqntd1pml} becomes  
\begin{equation}
\label{eqntd1pml1}
    \tilde{\mathcal{K}}_n \tilde{u}^s_n -\frac{1}{{2 \pi}} (\tilde{\mathcal{K}}^{0}1)_0 \tilde{u}^s_n= \tilde{\mathcal{S}}_n 
\partial_{\nu_c} \tilde{u}^s_n , \quad \text{on} \quad \gamma^a.
\end{equation}
Similarly, integrating equation \eqref{eq5} yields
\begin{equation}
    \label{eqntd1pml2}
    \tilde{\mathcal{K}}_{n,AB} \tilde{u}^s_n -\frac{1}{{2 \pi}} ({\mathcal{K}}_{AB}^{0}1 +{\mathcal{K}}_{+}^{0}1)_0 \tilde{u}^s_n= \tilde{\mathcal{S}}_{n,AB} 
\partial_{\nu_c} \tilde{u}^s_n , \quad \text{on} \quad \gamma_{AB}.
\end{equation}
where $\tilde{\mathcal{K}}_{n,AB}$, ${\mathcal{K}}_{n,AB}^0$ and $\tilde{\mathcal{S}}_{n,AB}$ are analogous to $\tilde{\mathcal{K}}_{n}$, ${\mathcal{K}}_{n}^0$ and $\tilde{\mathcal{S}}_{n}$, respectively, but with the integral domain $\gamma^a$ replaced by $\gamma_{AB}$. It is worth noting that ${\mathcal{K}}^{0}_{+}1 (\boldsymbol{r})$ is a known function and independent of $\theta$. So we let $R(r) $ denotes the function $ \frac{1}{{2\pi}}({\mathcal{K}}^{0}_{+}1 (\boldsymbol{r}))_0$. Consequently, the integral \eqref{eqntd1pml2} becomes
\begin{equation}
    \label{eq6}
    \tilde{\mathcal{K}}_{n,AB} \tilde{u}^s_n -(\frac{1}{{2\pi}}({\mathcal{K}}^{0}_{AB}1)_0 (r)+ R(r))\tilde{u}^s_n= \tilde{\mathcal{S}}_{n,AB} 
\partial_{\nu_c} \tilde{u}^s_n, \quad \text{on} \quad \gamma_{AB}.
\end{equation}
From \eqref{eq6}, we define the PML-transformed NtD map $\tilde{\mathcal{N}_n}= (\tilde{\mathcal{K}}_{n,AB} \tilde{u}^s_n -\frac{1}{{2\pi}}({\mathcal{K}}^{0}_{AB}1)_0 (r)- R (r) )^{-1}\tilde{\mathcal{S}}_{n,AB} $, mapping $\partial_{\nu_c} \tilde{u}^s_n$ to $ \tilde{u}^s_n$ on $\gamma_{AB}$.
Numerically discretizing the integral operators in \eqref{eq6} approximates the NtD map $\tilde{\mathcal{N}_n}$ on $\gamma_{AB}$. 

\section{Numerical implementation}
In this section, we consider the numerical discretizations of integral operators $\tilde{\mathcal{K}}_{n, AB}$, $\tilde{\mathcal{S}}_{n, AB}$ and $({\mathcal{K}}^{0}_{AB}1)_0$ on $\gamma_{AB}$. Let $r(s)=\{ (\rho(s),z(s)| 0  \leq s\leq  L_{AB} \}$ represent the parametric equation for $\gamma_{AB}$, where $s$ denotes the arclength parameter, and $L_{AB}$ is the length of $\gamma_{AB}$.
Since corners or edges may exist at the interface of the layered media, the derivatives of $\tilde{u}^s_n$ may exhibit corner singularities.
Following  \cite{colton2013}, we address these singularities through a graded mesh using a scaling function $s=\xi(t)$, $0 \leq t \leq 1 $. Specifically, for a smooth segment of the piecewise smooth curve $\gamma_{AB}$ corresponding to $s\in [s^0,s^1]$ and $t\in [t^0, t^1]$ such that $s^l=\xi(t^l)$, $l=1,2 $, where $s^0$ and $s^1$ correspond to two corners, we take 
\begin{equation}
    \label{eqsxi}
    s=\xi(t)=\frac{s^0 \xi_1^p+s^1 \xi_2^p}{\xi_1^p+\xi_2^p}, \quad t \in [t^0, t^1],
\end{equation}
where
\begin{equation}
    \xi_1=\left(\frac{1}{2}-\frac{1}{p}\right) \xi^3+\frac{\xi}{p}+\frac{1}{2}, \quad \xi_2=1-\xi_1, \quad \xi=\frac{2 t-\left(t^0+t^1\right)}{t^1-t^0}.
\end{equation}
Note that $p$ has been used in \eqref{eqsgm1} to define $\sigma_1$, and derivatives of $\xi(t)$ vanish at the corners up to order $p-1$. Points $t\in [0,1]$ are discretized by $\{ t_j=jh\}_{j=1}^{N}$ with grid size $h=1/(N-1)$, where $N$ is the number of discrete points, including corner points. The function $s = \xi(t)$ creates a graded mesh on $\gamma_{AB}$, creating discrete points that cluster near the corner points.

To simplify the notations, we use $r(t)$ to denote $r(\xi(t))$ and $r^{\prime}(t)$ to denote $\frac{d r}{d s}(\xi(t)) \xi^{\prime}(t)$ in the following.

\subsection{Approximating $\tilde{\mathcal{N}}_n$ on $\gamma_{AB}$}
To obtain an accurate discretization scheme, we first consider the exact evaluation of the PML-transformed modal Green’s function.
Recall that the kernels of operators  $\tilde{\mathcal{S}}( {\boldsymbol{r}},  {\boldsymbol{r}}^\prime)$ and $\tilde{\mathcal{K}} ( {\boldsymbol{r}},  {\boldsymbol{r}}^\prime)$ are denoted by $\tilde{S}( {\boldsymbol{r}},  {\boldsymbol{r}}^\prime)$ and $\tilde{K}( {\boldsymbol{r}},  {\boldsymbol{r}}^\prime)$, respectively. These kernels exhibit weakly singularity  at ${\boldsymbol{r}}={\boldsymbol{r}}^\prime$. To extract the most singular term of the kernels, we split $\tilde{S}({\boldsymbol{r}}, {\boldsymbol{r}}^\prime)$ and $\tilde{K}({\boldsymbol{r}}, {\boldsymbol{r}}^\prime)$ as 
\begin{align}
\label{eq:s1}
& \tilde{S}({\boldsymbol{r}}, {\boldsymbol{r}}^\prime)=Z\left(\tilde{\boldsymbol{r}},  \tilde{\boldsymbol{r}}^{\prime}\right)\left (H_1\left(\tilde{\boldsymbol{r}},  \tilde{\boldsymbol{r}}^{\prime}\right)+\bi H_2\left(\tilde{\boldsymbol{r}},  \tilde{\boldsymbol{r}}^{\prime} \right)\right),  \\
\label{eq:k1}
& \tilde{K}({\boldsymbol{r}}, {\boldsymbol{r}}^\prime)=D \left(\tilde{\boldsymbol{r}},  \tilde{\boldsymbol{r}}^{\prime}\right)\left (H_3\left(\tilde{\boldsymbol{r}},  \tilde{\boldsymbol{r}}^{\prime}\right)+{\bi H}_4\left(\tilde{\boldsymbol{r}},  \tilde{\boldsymbol{r}}^{\prime} \right)\right), 
\end{align}
where $\tilde{\boldsymbol{r}}=(\tilde{\rho}\cos\theta, \tilde{\rho}\sin\theta,\tilde{z})$, $\tilde{\boldsymbol{r}}'=(\tilde{\rho}'\cos\theta', \tilde{\rho}'\sin\theta',\tilde{z}')$,
\begin{align*}
& Z\left(\tilde{\boldsymbol{r}},  \tilde{\boldsymbol{r}}^{\prime}\right)=\frac{1}{4 \pi\left|\tilde{\boldsymbol{r}}-\tilde{\boldsymbol{r}}^{ \prime}\right|},  \\
& D \left(\tilde{\boldsymbol{r}},  \tilde{\boldsymbol{r}}^{\prime}\right)=-\frac{\nu_{c,\rho}^\prime \left(\tilde{\rho}^{\prime}-\tilde{\rho} \cos \left(\theta-\theta^{\prime}\right)\right)+\nu_{c,z}^{\prime}\left(\tilde{z}^{\prime}-\tilde{z}\right)}{4 \pi\left|\tilde{\boldsymbol{r}}-\tilde{\boldsymbol{r}}^{\prime}\right|^3},  \\
& H_1\left(\tilde{\boldsymbol{r}},  \tilde{\boldsymbol{r}}^{\prime}\right)=\cos \left(\omega_1\left|\tilde{\boldsymbol{r}}-\tilde{\boldsymbol{r}}^{\prime} \right|\right),  \\
& H_2\left(\tilde{\boldsymbol{r}},  \tilde{\boldsymbol{r}}^{\prime}\right)=\sin \left(\omega_1\left|\tilde{\boldsymbol{r}}-\tilde{\boldsymbol{r}}^{\prime} \right|\right),  \\
& H_3\left(\tilde{\boldsymbol{r}},  \tilde{\boldsymbol{r}}^{\prime}\right)=\cos \left(\omega_1\left|\tilde{\boldsymbol{r}}-\tilde{\boldsymbol{r}}^{\prime} \right|\right)+\omega_1\left|\tilde{\boldsymbol{r}}-\tilde{\boldsymbol{r}}^{\prime}\right| \sin \left(\omega_1\left|\tilde{\boldsymbol{r}}-\tilde{\boldsymbol{r}} ^{\prime}\right|\right),  \\
& H_4\left(\tilde{\boldsymbol{r}},  \tilde{\boldsymbol{r}}^{\prime}\right)=\sin \left(\omega_1\left|\tilde{\boldsymbol{r}}-\tilde{\boldsymbol{r}}^{\prime} \right|\right)-\omega_1\left|\tilde{\boldsymbol{r}}-\tilde{\boldsymbol{r}}^{\prime}\right| \cos \left(\omega_1\left|\tilde{\boldsymbol{r}}-\tilde{\boldsymbol{r}} ^{\prime}\right|\right).
\end{align*}
Integrating both sides of equations \eqref{eq:s1} and \eqref{eq:k1} gives 
\begin{align}
\label{eqsnspml}
    \tilde{S}_n(r,  r^{\prime})=& (ZH_1)_n(\tilde{r},\tilde{r}') + \bi(ZH_2)_n(\tilde{r},\tilde{r}'),\\
\label{eqknspml}
    \tilde{K}_n(r,  r^{\prime})=& (DH_3)_n(\tilde{r},\tilde{r}') + \bi(DH_4)_n(\tilde{r},\tilde{r}').
\end{align} 
Here, $\tilde{r}=(\tilde{\rho},\tilde{z})$ and $\tilde{r}'=(\tilde{\rho}',\tilde{z}')$ are the projections of $\tilde{\boldsymbol{r}}$ and $\tilde{\boldsymbol{r}}'$ in the $\rho$-$z$ half-plane, respectively. 

Clearly, the distance function $|\boldsymbol{\tilde{r}} -\boldsymbol{\tilde{r}}'|$ is complex-valued. 
When $\tilde{r}$ and $\tilde{r}'$ are well-separated, the kernels ${S}_n(\tilde{r},  \tilde{r}^{\prime})$ and ${K}_n(\tilde{r},  \tilde{r}^{\prime})$ are computed via FFT in the azimuthal direction. When $\tilde{r}$ and $\tilde{r}'$ are close, the function $H_1$, $H_3$, $ZH_2$, $DH_4$ are smooth while  $H_2$, $H_4$, $Z$, $ZH_1$, $D$, $DH_3$ are non-smooth. For a smooth function $f(\boldsymbol{r})$, its Fourier mode $f_n(\tilde{r})$ decays rapidly with $n$ and converges rapidly with discrete points number $N$ in the FFT. Therefore, the functions $(ZH_2)_n(\tilde{r},\tilde{r}')$ and $(DH_4)_n({r},{r}')$ are computed via FFT. Since the Fourier mode of $fg(\tilde{r})$ has a rapid convergence if at least one series of $f_n(\tilde{r})$ and $g_n(r)$ decays rapidly, the functions $(ZH_1)_n(\tilde{r},\tilde{r}')$ and $(DH_3)_n(\tilde{r},\tilde{r}')$ are computed via convolution.
For the convolution of the functions $(ZH_1)_n(\tilde{r},\tilde{r}')$ and $(DH_3)_n(\tilde{r},\tilde{r}')$, the required functions $Z_n(\tilde{r},\tilde{r}')$ and $D_n(\tilde{r},\tilde{r}')$ are evaluated by semi-analytical techniques and special functions. See  \cite{LAI2019152,helsing2014,olver2010} for further details on the approximation of the functions $Z_n(\tilde{r},\tilde{r}')$ and $D_n(\tilde{r},\tilde{r}')$.

Next, we discretize the integral equation along the generating curve $\gamma_{AB}$.
Since the modal Green's functions possess logarithmic singularities \cite{conway2010}, our discretization scheme is based on the hybrid Gauss-trapezoidal quadrature rules developed by Alpert \cite{Alpert1999}.
The operators $\tilde{\mathcal{S}}_{n,AB}$ and $\tilde{\mathcal{K}}_{n,AB}$ at $r=r(t_l)$, $l=1,...,N$, can be parameterized by
\begin{align}
    \label{eqsnab}
    \tilde{\mathcal{S}}_{n,AB}({\partial_{\nu_c} \tilde{u}^s_n})(r(t_l))&=\int_0^1 \tilde{S}_n(r(t_l),  r(t)) \varphi(t) dt, \\
    \tilde{\mathcal{K}}_{n,AB}({\tilde{u}^s_n})(r(t_l))&=\int_0^1 \tilde{K}_n(r(t_l),  r(t)) \tilde{u}^s_n dt,
\end{align}
where $\varphi(t) = {\partial_{\nu_c} \tilde{u}^s_n} |r'(t)|$. Taking $\tilde{\mathcal{S}}_{n,AB}({\partial_{\nu_c} \tilde{u}^s_n})(r(t_l))$ as an example, the integrand in \eqref{eqsnab} have logarithmic singularities at $t_l$. Therefore, we rewrite \eqref{eqsnab} as
\begin{equation}
    \label{eqsnab1}
    \tilde{\mathcal{S}}_{n,AB}({\partial_{\nu_c} \tilde{u}^s_n})(r(t_l))=\int_0^{t_l} \tilde{S}_n(r(t_l),  r(t)) \varphi(t) dt
    +\int_{t_l}^1 \tilde{S}_n(r(t_l),  r(t)) \varphi(t) dt.
\end{equation}
Let $N_1$ and $N_2$ denote the number of discrete points on intervals $[0,t_l]$ and $[t_l,1]$, respectively.

By applying Alpert’s hybrid Gauss-trapezoidal quadrature rule, the first part of the r.h.s of \eqref{eqsnab1} can be discretized as 
\begin{equation}
\label{eqsnt1}
\begin{aligned}
\int_0^{t_l} \tilde{S}_n(r(t_l),  r(t)) \varphi(t) dt\approx 
& \sum_{k=1}^{K_1} \gamma_k h \tilde{S}\left(t_l, \delta_k h\right) \varphi\left(\delta_k h\right) \\
& +\sum_{k=1}^{K_2} \gamma_k^\prime h \tilde{S}\left(t_l, t_l-\delta_k^\prime h\right) \varphi\left(t_l-\delta_k^\prime h\right) \\
& +\sum_{k=K_3}^{N_1-K_4} h \tilde{S}\left(t_l, t_k\right) \varphi\left(t_k\right),
\end{aligned}
\end{equation}
where the values of $K_1$, $K_2$, $K_3$, $K_4$, $\gamma_k$, $\gamma_k^\prime$, $\delta_k$ and $\delta_k^\prime$ are depend on the order of Alpert’s quadrature rule and can be precomputed. 
In cases where $N_1 < K_1+ K_2$, the third sum on the r.h.s of \eqref{eqsnt1} becomes irrational. In this situation, we rescaled the interval $[0,t_l]$ with a new grid size $h_1=1/(K_1+ K_2-1)$, ensuring that the number of uniform nodes on the interval $[0,t_l]$ is at least $K_1+ K_2$. The second sum on the r.h.s of equation \eqref{eqsnab1} can be discretized similarly.
For $l=1,...,N$, adopting the Lagrange interpolation formula to approximate $\varphi\left(\delta_k h\right)$, $\varphi\left(t_l-\delta_k^\prime h\right)$, $\varphi\left(t_l + \delta_k^\prime h\right)$ and $\varphi\left(1-\delta_k h\right)$ in terms of $\{ \varphi(t_j)\}_{j=1}^N$, we obtain
\begin{equation}
\left(\tilde{\mathcal{S}}_{n,AB}  \partial_{\nu_c} \tilde{u}^s_n\right)\left[\begin{array}{c}
r\left(t_1\right) \\
\vdots \\
r\left(t_N\right)
\end{array}\right] \approx  \tilde{\boldsymbol{S}}_n\left[\begin{array}{c}
\varphi\left(t_1\right) \\
\vdots \\
\varphi \left(t_N\right)
\end{array}\right],
\end{equation}
where $\tilde{\boldsymbol{S}}_n$ is an $N \times N$ matrix approximating the operator $\tilde{\mathcal{S}}_{n,AB}$. Similarly, an $N \times N$ matrix $\tilde{\boldsymbol{K}}_n$ that approximates the operator $\tilde{\mathcal{K}}_{n,AB}$ can be derived as
\begin{equation}
\left(\tilde{\mathcal{K}}_{n,AB}  \tilde{u}^s_n\right)\left[\begin{array}{c}
r\left(t_1\right) \\
\vdots \\
r\left(t_N\right)
\end{array}\right] \approx  \tilde{\boldsymbol{K}}_n\left[\begin{array}{c}
 \tilde{u}^s_n\left(r\left(t_1\right)\right) \\
\vdots \\
 \tilde{u}^s_n\left(r\left(t_N\right)\right)
\end{array}\right].
\end{equation}

The operator $(\frac{1}{{2\pi}}({\mathcal{K}}^{0}_{AB}1)_0 (r)+ R (r)$ in \eqref{eq6} can be approximated by an $N \times N$ diagonal matrix $\boldsymbol{D}$ with entries $(\frac{1}{{2\pi}}({\mathcal{K}}^{0}_{AB}1)_0 (r(t_j))+ R (r(t_j))$, $j=1,...,N$. Thus, equation \eqref{eq6} can be transformed in matrix equation form as
\begin{equation}
    (\tilde{\boldsymbol{K}}_n - \boldsymbol{D}) \tilde{\boldsymbol{u}}^s_n \approx \tilde{\boldsymbol{S}} \boldsymbol{\varphi}_n,
\end{equation}
where 
\begin{align}
    \tilde{{\bm u}}^s_n &= \left[u^s_n\left(r(t_1)\right),\cdots,u^s_n\left(r(t_N)\right) \right]^T,\\
    {\bm \varphi}_n &= \left[\varphi_n \left(r(t_1)\right),\cdots,\varphi_n\left(r(t_N)\right) \right]^T.
\end{align}
Consequently, one gets
\begin{equation}
        \tilde{\boldsymbol{u}}^s_n \approx \tilde{\boldsymbol{N}}_n \boldsymbol{\varphi}_n
        :=(\tilde{\boldsymbol{K}}_n - \boldsymbol{D})^{-1} \tilde{\boldsymbol{S}} \boldsymbol{\varphi}_n
\end{equation}
where the $N \times N$ matrix $\tilde{\boldsymbol{N}}_n$ approximates the NtD operator $\tilde{\mathcal{N}_n}$ mapping $\boldsymbol{\varphi}_n=|r'| \partial_{{\nu}_c} \tilde{\boldsymbol{u}}_n$ to $\tilde{\boldsymbol{u}}^s_n$ on $\gamma_{AB}$.

\subsection{Wave field evaluations}
Suppose now for each domain $\Omega_j$, we have obtained an $N \times N$ matrices $ \tilde{\boldsymbol{N}}^j_n$ approximating the NtD operator $\tilde{\mathcal{N}}_{n}^{j}$ mapping $|r'| \partial _{\nu_c} \tilde{\boldsymbol{u}}_{n}^{s,j}$ to $\tilde{\boldsymbol{u}}_{n}^{s,j}$ on $\gamma_{AB}$ for $j=1,2$.
Thus, we obtain the following relation
\begin{equation}
\label{eqntdnj}
    \tilde{\boldsymbol{u}}_{n}^{s,j} = \tilde{\boldsymbol{N}}_{n}^j \boldsymbol{\varphi}_{n}^j, \quad j=1,2, 
\end{equation}
where 
\begin{align}
    {\bm u}_{n}^{s,j} &= \left[u_{n}^{s,j}\left(r(t_1)\right),\cdots,u_{n}^{s,j}\left(r(t_N)\right) \right]^T,\\
    {\bm \varphi}_{n}^j &= \left[\varphi_{n}^j \left(r(t_1)\right),\cdots,\varphi_{n}^j\left(r(t_N)\right) \right]^T,
\end{align}
and $\varphi_{n}^j = |r'| \partial _{\nu_c} \tilde{\boldsymbol{u}}_{n}^{s,j}$.

Applying the transmission conditions \eqref{eqs1} and \eqref{eqs2}, the scattering wave field $\tilde{u}_{n}^{s,j}$ on $\gamma_{AB} $ satisfies 
\begin{align}
\tilde{{u}}_{n}^{s,1}-\tilde{{u}}_{n}^{s,2} & =\mathbf{b}_{n}^{1}, \\
{\bm \varphi}_{n}^{1}-{\bm \varphi}_{n}^{2} & =\mathbf{b}_{n}^{2},
\end{align}    
where
\begin{equation*}
\begin{aligned}
& \mathbf{b}_{n}^{1}=\left[-\left[\tilde{u}_{n}^{tot,0}\right]\left(r\left(t_1\right)\right), \ldots,-\left[\tilde{u}_{n}^{tot,0}\right]\left(r\left(t_N\right)\right)\right]^T, \\
& \mathbf{b}_{n}^{2}=\left[-\left|r^{\prime}\left(t_1\right)\right|\left[\partial_{\boldsymbol{\nu}_{\mathrm{c}}} \tilde{u}_{n}^{tot,0}\right]\left(r\left(t_1\right)\right), \ldots,-\left|r^{\prime}\left(t_N\right)\right|\left[\partial_{\boldsymbol{\nu}_{\mathrm{c}}} \tilde{u}_{n}^{tot,0}\right]\left(r\left(t_N\right)\right)\right]^T,
\end{aligned}    
\end{equation*}
and $\tilde{u}_{n}^{tot,0}$ denotes the $n$-th Fourier mode of  $\tilde{u}^{tot,0}$, which is the PML-transformed wave field of the reference wave field ${u}^{tot,0}$.
From equation \eqref{eqntdnj}, we obtain 
\begin{equation}
\left[\begin{array}{cc}
\tilde{\mathbf{N}}^1_{n} & -\tilde{\mathbf{N}}^2_{n} \\
 \mathbf{I} & -\mathbf{I}
\end{array}\right]\left[\begin{array}{l}
\boldsymbol{\varphi}_{n}^{1} \\
\boldsymbol{\varphi}_{n}^{2}
\end{array}\right]=\left[\begin{array}{l}
\mathbf{b}_{n}^{1} \\
\mathbf{b}_{n}^{2}
\end{array}\right],
\end{equation}
where $\mathbf{I}$ denotes the identity matrix and the solution
\begin{align}
\label{eqntdpsi1}
& \boldsymbol{\varphi}_{n}^{1}=\left(\tilde{\mathbf{N}}^1_{n}- \tilde{\mathbf{N}}^2_{n}\right)^{-1}\left( \mathbf{b}_{n}^{1}-\tilde{\mathbf{N}}^2_{n} \mathbf{b}_{n}^{2} \right), \\
\label{eqntdpsi2}
& \boldsymbol{\varphi}_{n}^{2}= \boldsymbol{\varphi}_{n}^{1}-{\mathbf{b}_{n}^{2}}.
\end{align}
Then, we obtain the scattering wave $\tilde{\boldsymbol{u}}_{n}^{s,j}=\tilde{\mathbf{N}}_{n}^j \boldsymbol{\varphi}_n^j$, for $j=1,2$, on $\gamma_{AB}$.
For $r\in \Omega_j$, using Green's representation formula, the scattering wave field is approximated as
\begin{equation}
   \label{eqfujs}
    \tilde{\boldsymbol{u}}_{n}^{s,j}(r) \approx \int_{\gamma_{A B}}\left\{\tilde{G}_j(r, r') \partial_{\boldsymbol{\nu}_c} \tilde{\boldsymbol{u}}_{n}^{s,j}(r')-\partial_{\boldsymbol{\nu}_c} \tilde{G}_j(r, r') \tilde{\boldsymbol{u}}_{n}^{s,j}(r')\right\} d s(r') .
\end{equation}
After parameterizing by the scaling function $s = w(t)$ from  \eqref{eqsxi}, the integral in \eqref{eqfujs} becomes sufficiently smooth. Applying the trapezoidal rule yields
\begin{equation}
\begin{aligned}
       \tilde{\boldsymbol{u}}_{n}^{s,j}(r) \approx &\frac{1}{N} \sum_{l=1}^N[\tilde{G}_j\left(r, r\left(t_l\right)\right)\left|r^{\prime}\left(t_l\right)\right| \partial_{\nu_c} \tilde{\boldsymbol{u}}_{n}^{s,j}\left(r\left(t_l\right)\right) \\
       &-\partial_{\nu_c} \tilde{G}_j^s\left(r, r\left(t_l\right)\right)\left|r^{\prime}\left(t_l\right)\right| \tilde{\boldsymbol{u}}_{n}^{s,j}\left(r\left(t_l\right)\right)] .
\end{aligned}
\end{equation}
The scattered field is then given by
${\boldsymbol{u}}_{n}^{s,j}=\tilde{\boldsymbol{u}}_{n}^{s,j}$. Consequently, in the physical domain outside the PML layer, the complete total wave field is obtained as
\begin{equation}
    {\boldsymbol{u}}^{{tot},j}={\boldsymbol{u}}^{{tot},0}+{\boldsymbol{u}}^{{s},j}, \quad \text{in} \quad \Omega_j,
\end{equation}
where $j=1,2$ and ${\boldsymbol{u}}^{{s},j}=\sum_{n=-\infty}^{\infty} {\boldsymbol{u}}_{n}^{{s},j}$.

\subsection{Instability issues}
In this subsection, we address several instability issues that arise during numerical computations and propose strategies to tackle them.
When evaluating the kernel functions $\tilde{S}_n(r, r')$ and $\tilde{K}_n(r, r')$, if $r$ is very close to $r'$ or close to the boundary, the distant function $|\tilde{\boldsymbol{r}}-\tilde{\boldsymbol{r}}^\prime|$ utilized in FFT may tend to zero, leading to numerical instabilities or potential divide-by-zero errors. To resolve these instabilities, we rewrite the distant function
\begin{align*}
    |\tilde{\boldsymbol{r}}-\tilde{\boldsymbol{r}}^\prime| &=\sqrt{\tilde{\rho}^2+\tilde{\rho}^{\prime 2}-2 \tilde{\rho} \tilde{\rho}^{\prime} + 4 \tilde{\rho} \tilde{\rho}^{\prime} \sin^2 \left(  \frac{\theta-\theta^{\prime}}{2}\right)+\left(\tilde{z}-\tilde{z}^{\prime}\right)^2} \\
    &=\sqrt{|\tilde{r}-\tilde{r}^\prime|^2 + 4 \tilde{\rho} \tilde{\rho}^{\prime} \sin^2 \left(  \frac{\theta-\theta^{\prime}}{2}\right)},
\end{align*}
Following the stabilization technique detailed in Section 4.3 of \cite{lulu2018}, we can compute $|\tilde{r}-\tilde{r}^\prime|$ accurately. Since the values of functions $|\tilde{r}-\tilde{r}^\prime|$ and $ 4 \tilde{\rho} \tilde{\rho}^{\prime} \sin^2 \left(  \frac{\theta-\theta^{\prime}}{2}\right)$ are both positive numbers, there are no significant rounding error. A similar approach was used to accurately evaluate $D \left(\tilde{\boldsymbol{r}}, \tilde{\boldsymbol{r}}^{\prime}\right)$.

Another instability issue arises during the convolution of the functions $(ZH_1)_n(\tilde{r},\tilde{r}')$ and $(DH_3)_n(\tilde{r},\tilde{r}')$. In the PML region, as $\theta \in [0,2\pi)$ varies, the positive imaginary part of the value for distance function $|\boldsymbol{\tilde{r}}-\boldsymbol{\tilde{r}}'|$ may become large so that the functions $H_1(\boldsymbol{\tilde{r}},\boldsymbol{\tilde{r}}')$ and $H_3(\boldsymbol{\tilde{r}},\boldsymbol{\tilde{r}}')$ grow exponentially, which leads to the series $H_{1n}(\tilde{r},\tilde{r}')$ and $H_{3n}(\tilde{r},\tilde{r}')$ converging at a slow rate. To overcome this problem, we introduce rapidly converging functions to replace the functions $H_1(\boldsymbol{\tilde{r}},\boldsymbol{\tilde{r}}')$ and $H_3(\boldsymbol{\tilde{r}},\boldsymbol{\tilde{r}}')$. Let the functions $H_1(\boldsymbol{\tilde{r}},\boldsymbol{\tilde{r}}')$ and $H_3(\boldsymbol{\tilde{r}},\boldsymbol{\tilde{r}}')$ be denoted by $H_{1}(q)$ and $H_{3}(q)$, respectively, where $q=\omega|\boldsymbol{\tilde{r}}-\boldsymbol{\tilde{r}}^{\prime}|$ is a complex number with a positive imaginary part. For simplicity, we consider the function $H_{1}(q)$ only. Applying the constant transformation
\begin{equation}
	H_{1}(q) = H_{1}(q) e^{-\bi q} e^{\bi q},
\end{equation}
the truncated Taylor expansion of the function $H_{1}(q) e^{-\bi q}$ is given by
\begin{equation}
    H_{1}(q) e^{-\bi q}=\sum_{k=0}^L a_k z^k+O\left(q^{L+1}\right),
\end{equation}
where $L$ denotes the truncation parameter. We define the function
\begin{equation}
\label{eqch1}
    c H_{1}(q, L)=e^{\bi q} \sum_{k=0}^L a_k z^k.
\end{equation}
Let the function $H_{1}(q)$ be replaced by the function $c H_{1}(q, L)$, which exhibits enhanced convergence properties. Accordingly, we replace $iH_{2}(q)$ with $e^{\bi q}-cH_{1}(q,L)$ to minimize the error raised by substitution. 

To illustrate the efficiency of this replacement, we carry out a specific experiment to validate the necessity of using the function $c H_{1}(q, L)$. As illustrated in Figure \ref{fg2}(b), the surface $\Gamma_{AB}$ is generated by rotating the line $\gamma_{AB}$ with endpoints ${\bm V}$ given by 
\begin{align*}
     {\bm V}&=\{(0,-1),(1,-1),(1,0),(4,0)\}.
\end{align*}
we set $\omega_1=\omega_2=2\pi$, $n=1$, $a_1=2$, the PML thickness $T=2$, and the PML parameter $S=2$. We have two choices to evaluate the kernel functions $\tilde{S}_n(r,r')$ and $\tilde{K}_n(r,r')$, where $q=\omega_1|\tilde{r}-\tilde{r}^\prime|$. First, using the original functions $H_{1}(q)$ and $H_{3}(q)$ to evaluate the kernel functions $\tilde{S}_n(r,r')$ and $\tilde{K}_n(r,r')$, respectively, we obtain an $N \times N$ matrix $ \tilde{\boldsymbol{N}}_n$ denoted by $\tilde{\boldsymbol{N}}^{\text{ori}}_n$. Second, we substitute the functions $cH_{1}(q,L)$ and $cH_{3}(q,L)$ for the original functions $H_{1}(q)$ and $H_{3}(q)$, respectively, to obtain the NtD matrix $\tilde{\boldsymbol{N}}^{\text{sub}}_n$.

\begin{figure}
    \centering
    \subfigure[]{\includegraphics[width=0.45\linewidth]{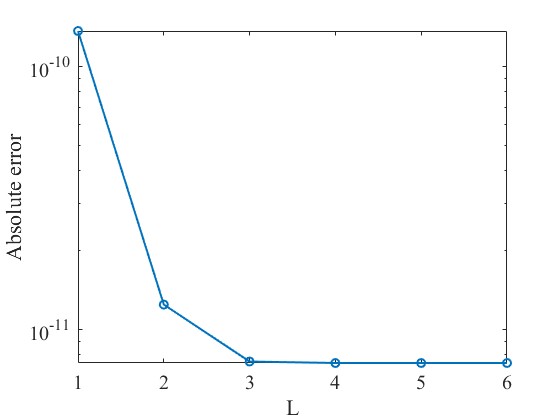}}
    \subfigure[]{\includegraphics[width=0.45\linewidth]{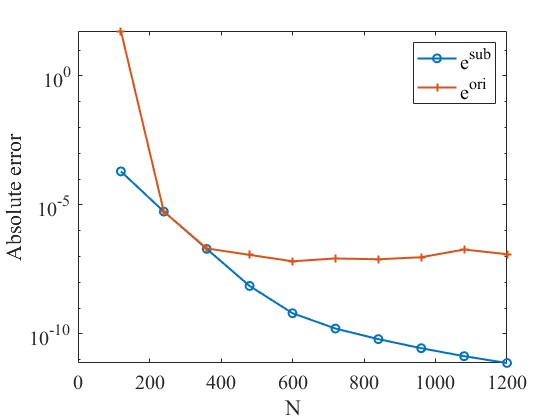}}
    \caption{(a) Accuracy of the NtD matrix ${\bm N}^{\text{sub}}_n$ for varying $L$. (b) Accuracy of the two NtD matrices ${\bm N}^{\text{ori}}_n$ and ${\bm N}^{\text{sub}}_n$, for $L=5$.}
    \label{picntde}
\end{figure}

To compare the accuracy of the two NtD matrices, we employ the following explicit solution 
\begin{equation}
\label{equp}
\tilde{u}({\bm r})= \tilde{S}\left({\boldsymbol{r}},  {\boldsymbol{r}}_{p}\right)= \frac{e^{\bi\omega_1|\tilde{{\bm r}}-\tilde{{\bm r}}_p|}}{4\pi|\tilde{{\bm r}}-\tilde{{\bm r}}_p|},   
\end{equation}
with ${\bm r}_p=(0.5,0,-1.3)$. Thus, the excited fields $\tilde{u}_n$ and the normal derivative $\partial_\nu \tilde{u}_n$ are given by
\begin{equation}
    \tilde{u}_n(r)=\tilde{S}_n(r,r_p) ,\quad
    \partial_\nu \tilde{u}_n(r)=\tilde{K}_n(r,r_p), \quad
    r\in \gamma_{AB}.
\end{equation}
The corresponding errors are defined as
\begin{equation}
    e^i=|| {\bm u}_n- \tilde{\bm N}^{i}_n{\bm \varphi}_n||_2,\quad i\in\{\text{ori},\text{sub}\}.
\end{equation} 
For $N=1200$, we plot the error function with varying $L$ in Figure \ref{picntde}(a). Accurate results can be produced even with a small $L$, and for convenience, we set $L=5$ in the following. The absolute errors of the two error functions against the number of discretization points $N$ are shown in Figure \ref{picntde}(b). Numerical results validate the necessity of using the function $cH_{1}(q,L)$. 

\section{Numerical examples}
In this section, we present several numerical experiments to illustrate the efficiency of our algorithm. The physical domain is defined as $\Omega_{phy} = \{(\rho, \theta, z) | \rho \leq a_1 , 0 \leq \theta<2\pi ,|z|\leq a_2  \}$. Since the integration domain $\gamma_{AB}$ is independent of $a_2$, we let $a_2$ approach infinity. Therefore, the PML domain $\Omega_{pml}=\{ (\rho, \theta, z) | a_1 \leq \rho \leq a_1+T, a_1 >0 , T>0 \}$.

In all examples, we take $p=6$ to define $\sigma_1$ and the scaling function $\xi(t)$, and apply Alpert's quadrature rule of the $10$-th order to discretize the integral equations. And we set $\omega_1=2\pi$, $\omega_2=2.4\pi$, $a_1=2$. In the following, we consider two kinds of incident waves: a plane wave incidence with the incident angle $\phi^{inc}=\pi/3$, $\theta^{inc}=0$, and a spherical wave incidence excited by the source $r^*=(0.3,0,0.5)$.

\begin{figure}
    \centering
    \subfigure[]{
    \includegraphics[width=0.31\linewidth]{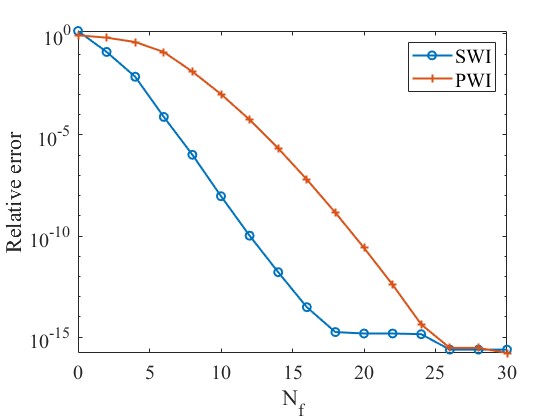}
    }
    \subfigure[]{
    \includegraphics[width=0.31\linewidth]{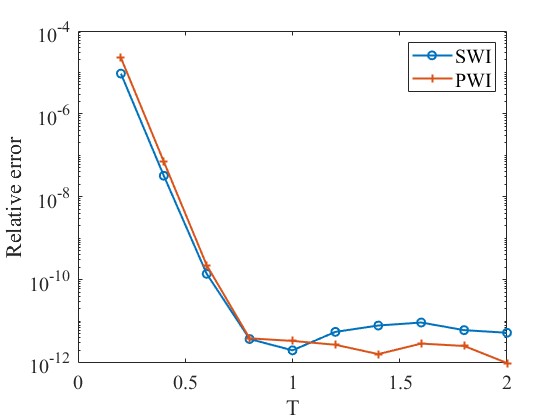}
    }
    \subfigure[]{
    \includegraphics[width=0.31\linewidth]{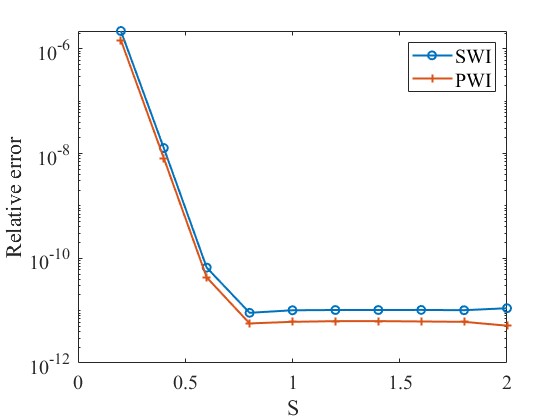}
    }
    \caption{Result for Example 1 with $N=800$. (a) Relative errors of $c\boldsymbol{u}^{s,j}$ with respect to $N_f$. (b) Relative errors of $u^{\text{tot}}$ against $T$ for $S=2$. (c) Relative errors of $u^{\text{tot}}$ against $S$ for $T=2$. SWI: Spherical wave incidence, PWI: Plane wave incidence.}
    \label{fg5}
\end{figure}
\begin{figure}
    \centering
    \subfigure[]{
    \includegraphics[width=0.45\linewidth]{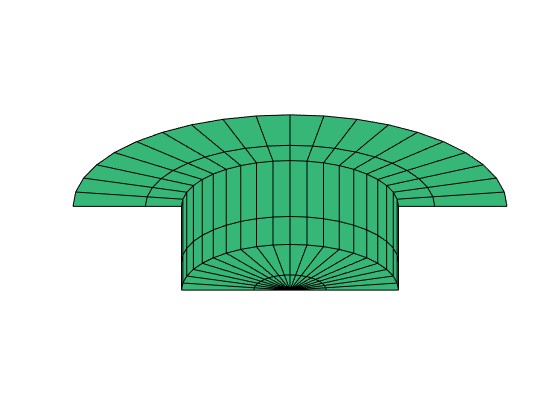}
    }
    \subfigure[]{
    \includegraphics[width=0.45\linewidth]{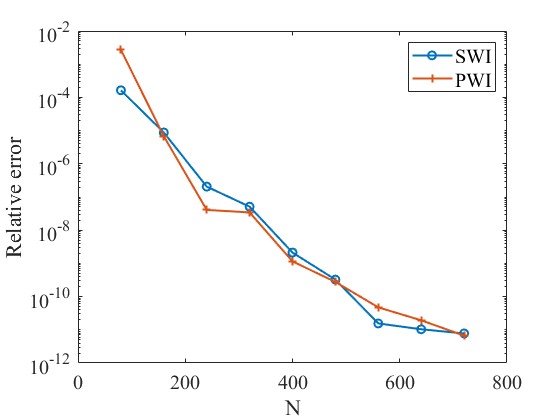}
    }
    \caption{Result for Example 1. (a)  Illustration of the perturbation. (b) Relative errors against $N$ for $T=2$ and $S=2$.   }
    \label{fg4}
\end{figure}
\begin{figure}
    \centering
    \subfigure[]{
    \includegraphics[width=0.45\linewidth]{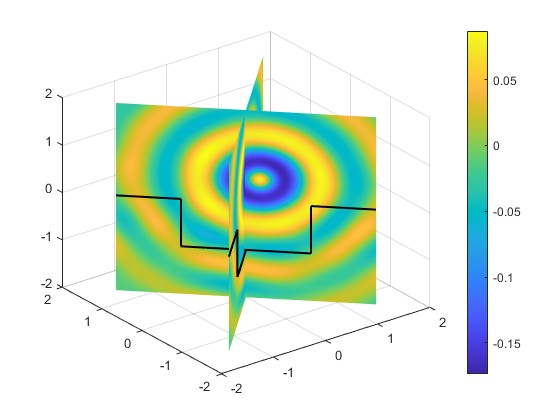}
    }
    \subfigure[]{
    \includegraphics[width=0.45\linewidth]{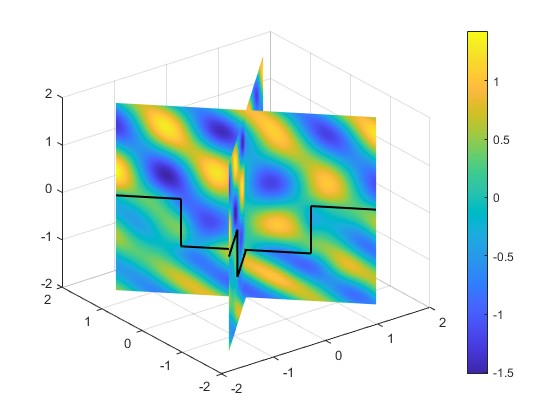}
    }

    \caption{Result for Example 1. (a) and (b) are the real parts of the scattered field for spherical and plane wave incidences, respectively.}
    \label{fg9}
\end{figure}

\textbf{Example 1.} In this example, we consider a perturbation with a rectangular cross-section, as shown in Figure \ref{fg4}(a). The depth and radius of the perturbation are both set to $1$. First, we investigate the convergence of the Fourier modes of the scattered field $\boldsymbol{u}^{s,j}=\sum_{n=-\infty}^{\infty} \boldsymbol{u}^{s,j}_n$. Let the truncated scattered field be denoted by $c\boldsymbol{u}^{s,j} =\sum_{n=-N_f}^{N_f} \boldsymbol{u}^{s,j}_n$, where $N_f$ denotes the truncation parameter. We set the number of discrete points $N=800$, the PML thickness $T=2$ and the PML parameter $S=2$. For the case of spherical wave incidence, we compute $c\boldsymbol{u}^{s,j}$ for different values of $N_f$ ranging from $0$ to $32$. Using the numerical solution $c\boldsymbol{u}^{s,j}$ with $N_f=32$ as a reference, we compute the relative errors for the less accurate numerical solutions with lower $N_f$. Similarly, for plane wave incidence, using $c\boldsymbol{u}^{s,j}$ with $N_f=32$ as a reference solution to compute the relative errors for lower $N_f$. The numerical results are shown in Figure \ref{fg5}(a). 
We can observe that the truncated scattered fields $c\boldsymbol{u}^{s,j}$ for two types of wave incidence converge rapidly as $N_f$ increases.
Therefore, we let $\boldsymbol{u}^{s,j} =\sum_{n=-N_f}^{N_f} \boldsymbol{u}^{s,j}_n$. In the following, we set the Fourier truncation parameters to $N_f=20$ and $N_f=30$ for the spherical and plane wave incidence cases, respectively.

Next, we investigate the convergence of the scattered field $\boldsymbol{u}^{tot}$ as the PML thickness $T$ and the PML parameter $S$ vary. Let $N=800$.
For both spherical and plane wave incidences, the reference solutions are the total wave field $\boldsymbol{u}^{tot}$ for $T=2.2$. These reference values are utilized to calculate the relative errors for the less accurate numerical solutions with lower $T$. Numerical results for both incident waves are depicted in Figure \ref{fg5}(b). Similarly, using the solution for $S=2.2$ as a reference, the relative errors for the less accurate numerical solutions with lower $S$ are shown in Figure \ref{fg5}(c). It can be observed that the relative errors decay exponentially at the beginning and then yield to the discretization error. Those relative errors demonstrate the efficiency of the proposed PML-based BIE solver for three-dimensional problems.

Then, we study the relative errors of $\boldsymbol{u}^{tot}$ as the number of discrete points $N$ varies. For both spherical and plane wave incidences, the reference solutions are the total wave field  $\boldsymbol{u}^{tot}$ for $N=800$. These reference values are used to calculate the relative errors for the less accurate numerical solutions with lower $N$. Numerical results for both incident waves are shown in Figure \ref{fg4}(b). Figure \ref{fg9} shows the spherical wave scattering with $41$ modes and the plane wave scattering with $61$ modes at four azimuthal directions, $\theta=\frac{\pi}{4},\frac{3\pi}{4},\frac{5\pi}{4},\frac{7\pi}{4}$.

\begin{figure}
    \centering
    \subfigure[]{
    \includegraphics[width=0.45\linewidth]{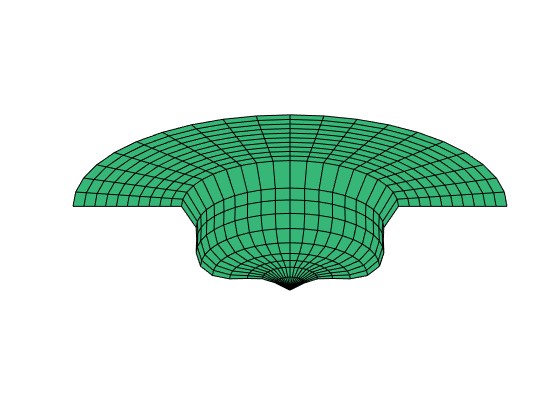}
    }
    \subfigure[]{
    \includegraphics[width=0.45\linewidth]{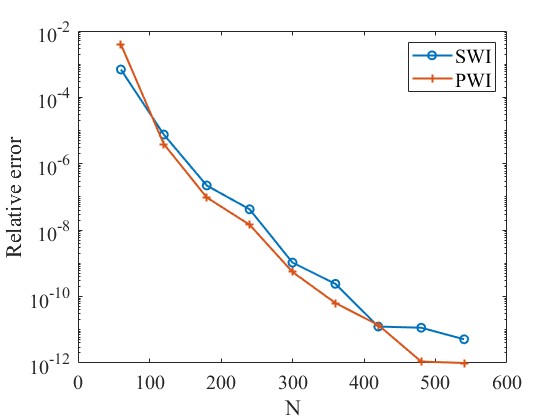}
    }

    \caption{Result for Example 2. (a)  Illustration of the perturbation. (b) Relative errors against $N$ for $T=2$ and $S=2$.}
    \label{fg6}
\end{figure}
\begin{figure}
    \centering
    \subfigure[]{
    \includegraphics[width=0.45\linewidth]{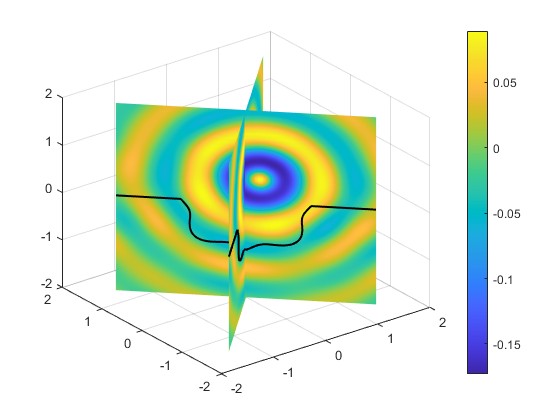}
    }
    \subfigure[]{
    \includegraphics[width=0.45\linewidth]{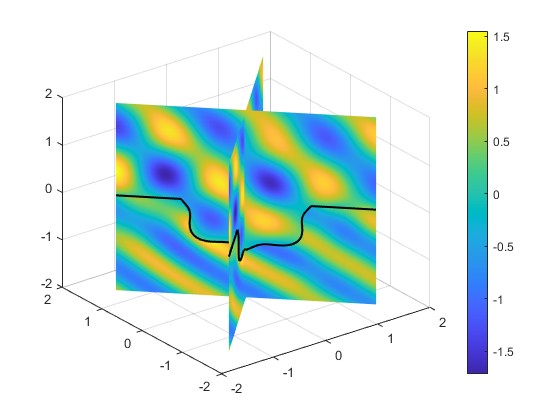}
    }

    \caption{Result for Example 2. (a) and (b) are the real parts of the scattered field for spherical and plane wave incidences, respectively.}
    \label{fg10}
\end{figure}
\textbf{Example 2.} In the second example, we consider an axisymmetric perturbation with the generating curve
\begin{equation}
\begin{aligned}
& \rho(s)=[1-0.1 \sin (3 \pi s)] \sin \left(\frac{\pi}{2} s\right), \\
& z(s)=-[1-0.1 \sin (3 \pi s)] \cos \left(\frac{\pi}{2} s\right),
\end{aligned}
\end{equation}
for $s\in [0,1]$, as depicted in \ref{fg6}(a). For two types of incident waves, we employ different values of $N_f$ like Example 1. Using $N=600$ points to discretize the generating cure $\gamma_{AB}$, we obtain a reference wave filed $\boldsymbol{u}^{tot}$. The relative errors for the less accurate numerical solutions with lower $N$ are calculated from this reference value. Numerical results for both incident waves are shown in Figure \ref{fg6}(b).  The scattered fields for spherical and plane wave incidence are given in Figure \ref{fg10}.

\begin{figure}
    \centering
    \subfigure[]{
    \includegraphics[width=0.4\linewidth]{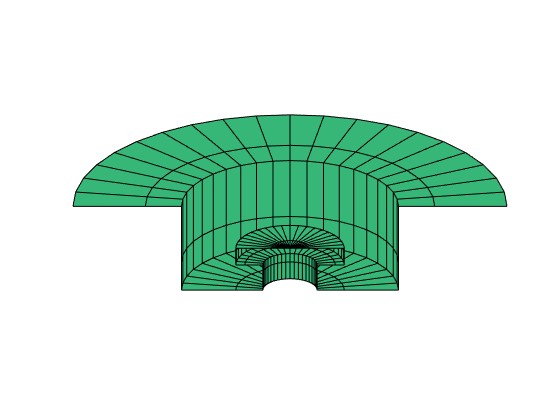}
    }
    \subfigure[]{
    \includegraphics[width=0.4\linewidth]{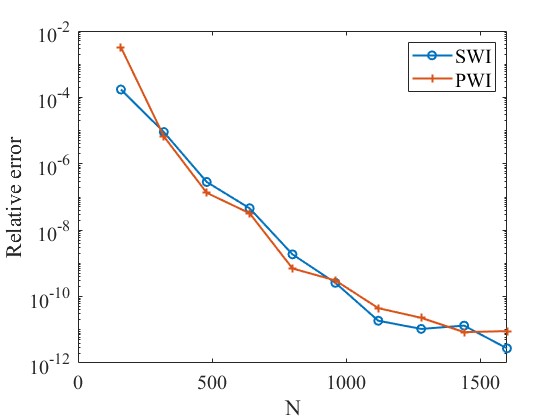}
    }
 
    \caption{Result for Example 3. (a)  Illustration of the perturbation. (b) Relative errors against $N$ for $T=2$ and $S=2$.}
    \label{fg7}
\end{figure}
\begin{figure}
    \centering
    \subfigure[]{
    \includegraphics[width=0.45\linewidth]{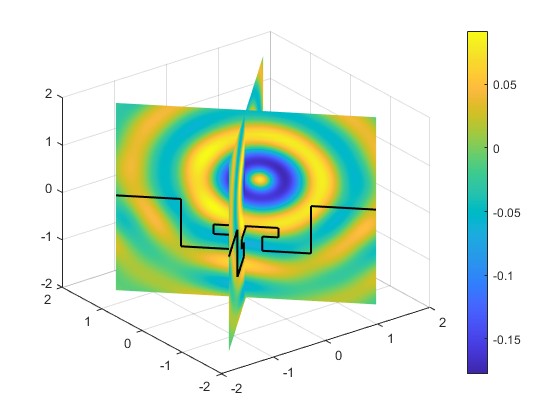}
    }
    \subfigure[]{
    \includegraphics[width=0.45\linewidth]{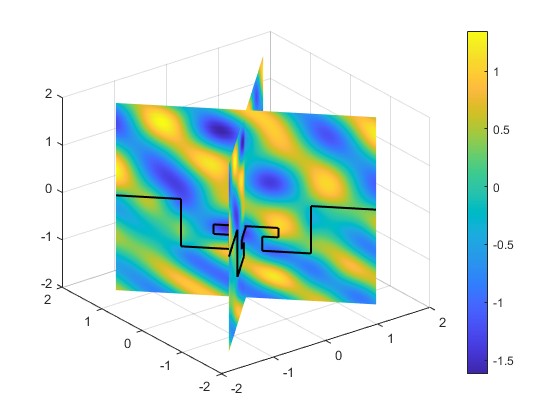}
    }

    \caption{Result for Example 3. (a) and (b) are the real parts of the scattered field for spherical and plane wave incidences, respectively.}
    \label{fg11}
\end{figure}
\textbf{Example 3.} In the third example, we consider a perturbation whose generating curve has vertices given by 
\begin{equation}
V=\{(0,-0.5),(0.5,-0.5),(0.5,-0.7),(0.25,-0.7),(0.25,-1),(1,-1),(1,0)\},
\end{equation}
as shown in Figure \ref{fg7}(a). Each segment of the piecewise smooth is rescaled by the function $\xi(t)$ to resolve the corner singularities. Using $N=1600$ points to discretize the generating cure $\gamma_{AB}$, we obtain a reference wave filed $\boldsymbol{u}^{tot}$. The relative errors for the less accurate numerical solutions with lower $N$ are calculated by this reference value. Numerical results for both incident waves are shown in Figure \ref{fg7}(b). Figure \ref{fg11} gives the scattered field for spherical and plane wave incidence.

\bibliographystyle{plain}
\bibliography{sample}

\begin{thebibliography}{10}

\bibitem{monro2008}
Monro~J. A.
\newblock {\em A super-algebraically convergent, windowing-based approach to the evaluation of scattering from periodic rough surfaces}.
\newblock California Institute of Technology, 2008.

\bibitem{Alpert1999}
B.~K. Alpert.
\newblock Hybrid gauss-trapezoidal quadrature rules.
\newblock {\em SIAM Journal on Scientific Computing}, 20(5):1551--1584, 1999.

\bibitem{Arens05}
T.~Arens and T.~Hohage.
\newblock On radiation conditions for rough surface scattering problems.
\newblock {\em IMA Journal of Applied Mathematics}, 70(6):839--847, 12 2005.

\bibitem{bruno2012regularized}
O.~Bruno, T.~Elling, and C.~Turc.
\newblock Regularized integral equations and fast high-order solvers for sound-hard acoustic scattering problems.
\newblock {\em International Journal for Numerical Methods in Engineering}, 91(10):1045--1072, 2012.

\bibitem{bruno2016}
O.~P Bruno, M.~Lyon, C.~P{\'e}rez-Arancibia, and C.~Turc.
\newblock Windowed green function method for layered-media scattering.
\newblock {\em SIAM Journal on Applied Mathematics}, 76(5):1871--1898, 2016.

\bibitem{bruno2021}
O.~P. Bruno and T.~Yin.
\newblock A windowed green function method for elastic scattering problems on a half-space.
\newblock {\em Computer Methods in Applied Mechanics and Engineering}, 376:113651, 2021.

\bibitem{cai2002}
W.~Cai.
\newblock Algorithmic issues for electromagnetic scattering in layered media: Green's functions, current basis, and fast solver.
\newblock {\em Advances in Computational Mathematics}, 16:157--174, 2002.

\bibitem{chandler98}
S.~N. Chandler-Wilde and B.~Zhang.
\newblock A uniqueness result for scattering by infinite rough surfaces.
\newblock {\em SIAM Journal on Applied Mathematics}, 58(6):1774--1790, 1998.

\bibitem{chen2010}
Z.~Chen and W.~Zheng.
\newblock Convergence of the uniaxial perfectly matched layer method for time-harmonic scattering problems in two-layered media.
\newblock {\em SIAM Journal on Numerical Analysis}, 48(6):2158--2185, 2010.

\bibitem{chew1999waves}
W.~C. Chew.
\newblock {\em Waves and fields in inhomogenous media}, volume~16.
\newblock John Wiley \& Sons, 1999.

\bibitem{Chew1994}
W.~C. Chew and W.~H. Weedon.
\newblock A 3d perfectly matched medium from modified maxwell's equations with stretched coordinates.
\newblock {\em Microwave and Optical Technology Letters}, 7(13):599--604, 1994.

\bibitem{colton2013}
D.~Colton and R.~Kress.
\newblock {\em Inverse Acoustic and Electromagnetic Scattering Theory}.
\newblock Applied Mathematical Sciences. Springer Berlin Heidelberg, 2013.

\bibitem{conway2010}
J.~T. Conway and H.~S. Cohl.
\newblock Exact fourier expansion in cylindrical coordinates for the three-dimensional helmholtz green function.
\newblock {\em Zeitschrift f{\"u}r angewandte Mathematik und Physik}, 61:425--443, 2010.

\bibitem{gao2022}
Y.~Gao and W.~Lu.
\newblock Wave scattering in layered orthotropic media i: a stable pml and a high-accuracy boundary integral equation method.
\newblock {\em SIAM Journal on Scientific Computing}, 44(4):B861--B884, 2022.

\bibitem{helsing2014}
J.~Helsing and A.~Karlsson.
\newblock An explicit kernel-split panel-based nyström scheme for integral equations on axially symmetric surfaces.
\newblock {\em Journal of Computational Physics}, 272:686--703, 09 2014.

\bibitem{helsing2017}
J.~Helsing and A.~Karlsson.
\newblock Resonances in axially symmetric dielectric objects.
\newblock {\em IEEE Transactions on Microwave Theory and Techniques}, 65(7):2214--2227, July 2017.

\bibitem{LAI2019152}
J.~Lai and M.~O'Neil.
\newblock An fft-accelerated direct solver for electromagnetic scattering from penetrable axisymmetric objects.
\newblock {\em Journal of Computational Physics}, 390:152--174, 2019.

\bibitem{Lassas2001}
M.~Lassas and E.~Somersalo.
\newblock Analysis of the pml equations in general convex geometry.
\newblock {\em Proceedings of the Royal Society of Edinburgh: Section A Mathematics}, 131:1183 -- 1207, 10 2001.

\bibitem{lu2021}
W.~Lu.
\newblock Mathematical analysis of wave radiation by a step-like surface.
\newblock {\em SIAM Journal on Applied Mathematics}, 81(2):666--693, 2021.

\bibitem{lulu2018}
W.~Lu, Y.~Y. Lu, and J.~Qian.
\newblock Perfectly matched layer boundary integral equation method for wave scattering in a layered medium.
\newblock {\em SIAM Journal on Applied Mathematics}, 78(1):246--265, 2018.

\bibitem{luxu2023}
W.~Lu, L.~Xu, T.~Yin, and L.~Zhang.
\newblock A highly accurate perfectly-matched-layer boundary integral equation solver for acoustic layered-medium problems.
\newblock {\em SIAM Journal on Scientific Computing}, 45(4):B523--B543, 2023.

\bibitem{meier2001}
A.~Meier and S.~N. Chandler-Wilde.
\newblock On the stability and convergence of the finite section method for integral equation formulations of rough surface scattering.
\newblock {\em Mathematical Methods in the Applied Sciences}, 24(4):209--232, 2001.

\bibitem{miret2013}
D.~Miret, G.~Soriano, and M.~Saillard.
\newblock Rigorous simulations of microwave scattering from finite conductivity two-dimensional sea surfaces at low grazing angles.
\newblock {\em IEEE Transactions on Geoscience and Remote Sensing}, 52(6):3150--3158, 2013.

\bibitem{monk2003finite}
P.~Monk.
\newblock {\em Finite element methods for Maxwell's equations}.
\newblock Oxford university press, 2003.

\bibitem{olver2010}
F.~Olver, D.~Lozier, R.~Boisvert, and C.~Clark.
\newblock Nist handbook of mathematical functions.
\newblock {\em US Department of Commerce, National Institute of Standards and Technology}, 2010.

\bibitem{ONEIL2018263}
M.~O'Neil and A.~J. Cerfon.
\newblock An integral equation-based numerical solver for taylor states in toroidal geometries.
\newblock {\em Journal of Computational Physics}, 359:263--282, 2018.

\bibitem{Paulus2000}
M.~Paulus, P.~Gay-Balmaz, and O.~J. Martin.
\newblock Accurate and efficient computation of the green’s tensor for stratified media.
\newblock {\em Physical Review E}, 62(4):5797, 2000.

\bibitem{perez2014}
C.~P{\'e}rez-Arancibia and O.~P. Bruno.
\newblock High-order integral equation methods for problems of scattering by bumps and cavities on half-planes.
\newblock {\em JOSA A}, 31(8):1738--1746, 2014.

\bibitem{saillard2011}
M.~Saillard and G.~Soriano.
\newblock Rough surface scattering at low-grazing incidence: A dedicated model.
\newblock {\em Radio Science}, 46(05):1--8, 2011.

\bibitem{Sommerfeld1909}
A.~Sommerfeld.
\newblock {\"U}ber die ausbreitung der wellen inder drahtlosen telegraphie.
\newblock {\em Annalen der Physik}, 333:665--736, 1909.

\bibitem{spiga2008}
P.~Spiga, G.~Soriano, and M.~Saillard.
\newblock Scattering of electromagnetic waves from rough surfaces: A boundary integral method for low-grazing angles.
\newblock {\em IEEE Transactions on Antennas and Propagation}, 56(7):2043--2050, 2008.

\bibitem{wanglu2024}
H.~Wang and W.~Lu.
\newblock A high-accuracy mode solver for acoustic scattering by a periodic array of axially symmetric obstacles.
\newblock {\em Journal of Scientific Computing}, 101(1):23, 2024.

\bibitem{yu2022}
X.~Yu, G.~Hu, W.~Lu, and A.~Rathsfeld.
\newblock Pml and high-accuracy boundary integral equation solver for wave scattering by a locally defected periodic surface.
\newblock {\em SIAM Journal on Numerical Analysis}, 60(5):2592--2625, 2022.

\bibitem{zhang2021accurate}
L.~Zhang, L.~Xu, and T.~Yin.
\newblock An accurate hypersingular boundary integral equation method for dynamic poroelasticity in two dimensions.
\newblock {\em SIAM Journal on Scientific Computing}, 43(3):B784--B810, 2021.

\bibitem{zhao2005}
Z.~Zhao, L.~Li, J.~Smith, and L.~Carin.
\newblock Analysis of scattering from very large three-dimensional rough surfaces using mlfmm and ray-based analyses.
\newblock {\em IEEE Antennas and Propagation Magazine}, 47(3):20--30, 2005.

\end{thebibliography}
\end{document}